\documentclass[twoside, 10pt]{article}

\usepackage{amsmath}
\usepackage{paralist}
\usepackage{amssymb,mathrsfs,psfrag,graphicx}
\usepackage{hyperref}
\def\e{\epsilon}

\def\t{\theta}
\def\T{\Theta}
\def\k{\kappa}

\def\a{\alpha}

\def\g{\gamma}
\def\d{\delta}

\def\s{\sigma}
\def\z{\zeta}

\def\di{\displaystyle}

\newtheorem{theorem}{Theorem}[section]

\newtheorem{lemma}[theorem]{Lemma}
\newtheorem{proposition}{Proposition}[section]

\newtheorem{remark}{Remark}

\begin{document}

\title{\bf   Stability of Superposition of Viscous Contact Wave and Rarefaction Waves
for Compressible Navier-Stokes System with Degenerate Heat-Conductivity and Large-Data} \vskip 0.5cm

\author{Manyu Liu$^a$
	\thanks{
		Email addresses: liumy@email.ncu.edu.cn(M. Y. Liu) }\\[3mm]  a. School of Mathematics and Computer Science ,\\
	University, Nanchang 330031, P. R. China }

\date{}
\maketitle
\begin{abstract}

This paper is concerned with the large-time behavior of solutions for the
compressible Navier-Stokes system with degenerate heat-conductivity describing the one-dimensional
motion of a viscous heat-conducting perfect polytropic gas. We proved that for the one-dimensional compressible system with temperature-dependent heat conductivity, the combination of viscous contact wave with rarefaction waves
for the non-isentropic  polytropic gas is asymptotically stable under large initial perturbation, provided the strength of the combination waves is suitably small.


\

 \noindent{\it Key
words and phrases:} Compressible Navier-Stokes system, Viscous Contact Discontinuity, 
 Degenerate Heat-Conductivity, Large Initial Perturbation,  Stability
\end{abstract}
\section{Introduction}
\renewcommand{\theequation}{\arabic{section}.\arabic{equation}}
\setcounter{equation}{0}

The one-dimensional compressible Navier-Stokes system describing the motion of a viscous heat-conducting perfect polytropic gas has the form in Lagrangian coordinates(see\cite{Batchelor,Serrin}):
\begin{equation}\label{ns}
\left\{
\begin{array}{ll}
\di v_t-u_x=0,\\
\di u_t+p_x=\mu\left(\frac{u_x}{v}\right)_x,\\
\di \left(e+\frac{u^2}{2}\right)_t+\left(p u\right)_x=\left(\kappa \frac{\t_x}{v}+\mu\frac{uu_x}{v}\right)_x
\end{array}
\right.
\end{equation}
for $x\in\mathbb{R}=(-\infty,+\infty)$, $t>0$, where $v(x,t)>0$, $u(x,t)$, $\t(x,t)>0$, $e(x,t)>0$ and $p(x,t)$
are the specific volume, fluid velocity, absolute temperature, internal energy
and pressure, respectively, while the positive constants $\mu$ and $\k$ denote
the viscosity and heat conduction coefficients, respectively. Here we study the ideal
polytropic fluids so that $p$ and $e$ are given by the state equations
\begin{equation}
	\label{p-e}
\di p=\frac{R\t}{v}=Av^{-\gamma}\exp\left(\frac{\gamma-1}{R}s\right),\quad e=c_{\nu}\t+\mathrm{const},
\end{equation}
where $s$ is the entropy, $\gamma>1$ is the adiabatic exponent,  $c_{\nu}=\frac{R}{\gamma-1}$ is
the specific heat, $A$ and $R$ are both positive constants. For $\mu$ and $\kappa$, we consider the case where $\mu$ and $\kappa$ are proportional to (possible different) powers of $\theta$:
\begin{equation}
	\di\mu = \tilde{\mu}\theta^\alpha,\quad \kappa = \tilde {\kappa}\theta^\beta,
	\label{mu-k}
\end{equation}
where $\tilde{\mu},\tilde{\kappa}>0$ and $\alpha, \beta \geq 0$ are constants.

Then we impose the Cauchy problem to the system \eqref{ns} supplement with the following initial and far field conditions:
\begin{equation}\label{initial}
\left\{
\begin{array}{ll}
\di (v,u,\t)(x,0)=(v_0,u_0,\t_0)(x), &\di x\in \mathbb{R},\\
\di (v,u,\t)(\pm\infty,t)=(v_{\pm},u_{\pm},\t_{\pm}), &\di t>0,
\end{array}
\right.
\end{equation}
where $v_{\pm}(>0)$, $u_{\pm}$ and $\t_{\pm}(>0)$ are given constants, and we assume
$\inf_{\mathbb{R}}v_0>0$, $\inf_{\mathbb{R}}\t_0>0$, and $(v_0,u_0,\t_0)(\pm\infty)=(v_{\pm},u_{\pm},\t_{\pm})$
as compatibility conditions. 

When the far field states are the same( $v_+=v_-$, $u_+=u_-$, $\t_+=\t_-$), a large number of literatures have been devoted to  some researches on the global existence of solutions to the system \eqref{ns} has advanced substantially since 1977, see \cite{jiang,jiang-2,kazhikhov,kazhi-she,L-L} and the reference therein. In particular, Jiang \cite{jiang, jiang-2} first obtained uniformly bounded for specific volume, while Li and Liang \cite{L-L} subsequently extended this framework by proving the uniformly estimates for temperature, thus completing the large-time behavior of solutions. 

It is interesting to study the global solutions in time of system \eqref{ns}
 and their large time behavior for different end states, which is more complicated than the same states. 
Notably, if the dissipation effects are neglected
($\mu=\k=0$), the system \eqref{ns} is reduced to the compressible Euler equations,
which is one of the most important hyperbolic system of conservation laws. It is well known that
the Euler system
 contains three basic wave patterns (see \cite{Smoller}), two nonlinear waves:
shock and rarefaction wave, and a linearly degenerate wave: contact discontinuity. When we consider the Riemann initial data
,the solutions consist of the above three wave patterns
and their superpositions, called by Riemann solutions, and govern both the local and large time
asymptotic behavior of general solutions of the Euler system. It is of great importance and interest
to study the large-time behavior of the viscous version of these basic wave patterns and their superpositions
to the compressible Navier-Stokes system \eqref{ns}.

Extensive literature has been devoted to the stability analysis of viscous wave patterns in the system \eqref{ns}. These works encompass the studies of \cite{Huang-matsumura,KM,MN-85} for shock waves, \cite{LX-88,MN-86,MN-92,nishi-yang-zhao} for rarefaction waves, and the reference therein. For the subtler case where the Riemann solution consists of contact discontinuity (referred to as viscous contact wave), Huang et al. in \cite{Huang-Matsumura-Xin,Huang-Xin-Yang} obtained the stability of viscous contact wave and convergence rate, and further related studies about stability of viscous contact wave can be found in \cite{huang-ma-shi,Huang-Yang,huang-zhao}. Subsequently, Huang-Li-Matsumura in \cite{Huang-Li-matsumura} proved the combination wave of a viscous contact wave with rarefaction waves is asymptotically stable in the  one-dimensional compressible Navier-Stokes system. Building on this, Huang-Wang \cite{H-W} succeeded in removing the non-physical condition of the adiabatic exponent $\gamma$, which means these basic wave patterns and their linear superpositions are stable even for large perturbation.

It is noted that the above mentioned literature has only considered cases where $\mu$ and $\kappa$ are constants, whereas according to the Chapman-Enskog expansion for the first level of approximation in kinetic theory,  the viscosity $\mu$ and heat conductivity $\kappa$ are functions of temperature alone \cite{Cercignani,Chapman}. When concerning the case where $\alpha,\beta>0$, the research becomes more challenging due to the
potential degeneracy and strong nonlinearity on viscosity and heat diffusion. For the existence of solutions, Jenssen-Karper \cite{Jenssen} showed the global existence of weak solutions when $\alpha = 0$
and $\beta\in (0, 3/2)$ on the bounded domain. Later on, for $\alpha=0$ and $\beta>0$, Pan-Zhang\cite{Pan-Zhang} (also cf.\cite{Duan-Guo-Zhu,Huang-Shi-Sun})obtained the
global strong solution. When the domain is unbounded, Li-Shu-Xu \cite{Li-Shu-Xu} proved the existence of global strong solutions under the initial conditions. For the large-time behavior of strong solutions with $\alpha =0$, $\beta>0$, when the domain is bounded,  Huang-Shi \cite{Huang-Shi} showed that the solution is nonlinearly exponentially stable as time tends to infinity. Recently, for the case of unbounded
domains, Li-Xu \cite{Li-Xu} made studies  on  the large-time behavior of the global strong solutions. However, the above literature has focused only on cases where the far field states are the same, whereas the cases of different far field states with temperature-dependent heat conductivity are still open.

The main aim of this paper is to prove the asymptotic stability of linear combination wave of viscous contact wave and the rarefaction waves of one-dimensional compressible system \eqref{ns}-\eqref{p-e} with temperature-dependent heat conductivity, provided the strength of combination waves is suitably small.

\emph{Notations.} Throughout this paper, generic positive
constants are denoted by $c$ and $C$ without confusion.

$L^{p}(\Omega)$ denotes the usual Lebesgue space on $\Omega
\subset\mathbb{R}=(-\infty,\infty)$ with its norm given by
\begin{equation*}
	\begin{array}{ll}
		\di\parallel f\parallel_{L^{p}(\Omega)}:=\left(\int_{\Omega}|f(x)|^{p}dx\right)^{\frac{1}{p}},
		\quad 1\leq p<\infty.\\[5mm]
		\di\parallel
		f\parallel_{L^{\infty}(\Omega)}:=\mbox{ess.sup}_{\Omega} |f(x)|,\quad\quad p=\infty.
	\end{array}
\end{equation*}

$H^{k}(\Omega)$ denotes the $k^{th}$ order Sobolev space with its
norm
$$
\|f\|_{H^{k}(\Omega)}:=\left(\sum ^{k}_{j=0}
\parallel \partial^{j}_{x}f\parallel^{2}(\Omega)\right)^{\frac{1}{2}}, \quad \mathrm{when} \parallel
\cdot
\parallel=\parallel \cdot
\parallel_{L^{2}(\Omega)}.
$$

The domain $\Omega$ will be often abbreviated without confusion.

\

\section{Main Result}

Prior to state the main results, we first recall the viscous contact wave $(V,U,\Theta)$
for the compressible Navier-Stokes system \eqref{ns} defined in \cite{Huang-Xin-Yang}. For the Euler system
\begin{equation}\label{euler}
	\left\{
	\begin{array}{ll}
		\di v_t-u_x=0,\\
		\di u_t+p_x=0,\\
		\di \left(e+\frac{u^2}{2}\right)_t+\left(p u\right)_x=0,
	\end{array}
	\right.
\end{equation} with Riemann
problem
\begin{equation}\label{Riemann}
	(v,u,\t)(x,0)=\left\{
	\begin{array}{ll}
		\di (v_-,u_-,\t_-),&\di x<0,\\
		\di (v_+,u_+,\t_+),&\di x>0,
	\end{array}
	\right.
\end{equation}it is known that the contact discontinuity solution takes the form
\begin{equation}
(\tilde V,\tilde U,\tilde\Theta)(x,t)=\left\{
\begin{array}{ll}
\di (v_-,u_-,\t_-),&\di x<0,t>0,\\
\di (v_+,u_+,\t_+),&\di x>0,t>0,
\end{array}
\right.
\end{equation}
provide that
\begin{equation}\label{RH}
\di u_-=u_+,\quad p_-\triangleq\frac{R\theta_-}{v_-}=p_+\triangleq\frac{R\t_+}{v_+}.
\end{equation}
We assume that $u_-=u_+=0$ without loss of generality.
Due to the dissipation effect, the corresponding wave $(V,U,\T)$ to the contact discontinuity $(\tilde V,\tilde U,\tilde \Theta)$ is smoothed and behaves as a diffusion wave, which is called a "viscous contact wave". The viscous contact wave $(V,U,\Theta)$ can be constructed
as follows. Since the pressure for the profile $(V,U,\Theta)$ is expected to be
constant asymptotically, we set
$$
\frac{R\Theta}{V}=p_+,
$$
which indicates the leading part of the energy equation $\eqref{ns}_3$ is
\begin{equation}\label{Theta t}
\di c_{\nu}\Theta_t+p_+U_x=\k_1\left(\frac{\Theta_x}{V}\right)_x,
\end{equation}
where $\k_1 = \tilde{\k}\T^\beta$. The equation \eqref{Theta t} and $\eqref{ns}_1$ lead to a nonlinear diffusion equation,
\begin{equation}\label{1.8}
\di \Theta_t=a\left(\Theta_x/\Theta\right)_x,\quad
\Theta(\pm\infty,t)=\theta_{\pm},\quad a=\frac{\kappa_1(\gamma-1) p_+}{\gamma R^2}>0,
\end{equation}
which has a unique self-similar solution $\Theta(x,t)=\Theta(\xi)$, $\xi=\frac{x}{\sqrt{1+t}}$
due to \cite{xiao-liu}. Furthermore, on the one hand, $\Theta(\xi)$ is a monotone function,
increasing if $\theta_+>\theta_-$ and decreasing if $\theta_+<\theta_-$; On the other hand, there exists some positive constant $\delta$, such that for $\d=|\t_+ -\t_-|$, $\Theta$ satisfies
\begin{equation}\label{1.9}
\di (1+t)|\Theta_{xx}|+(1+t)^{\frac{1}{2}}|\Theta_x|+|\Theta-\theta_{\pm}|
=O(1)\delta e^{-\frac{c_1x^2}{1+t}} \quad \mathrm{as}~ |x|\rightarrow\infty,
\end{equation}
where $c_1$ is positive constant depending only on $\t_{\pm}$. Once $\Theta$ is determined,
the contact wave profile $(V,U,\Theta)(x,t)$ is then defined as follows:
\begin{equation}\label{contact}
\di V=\frac{R\T}{p_+},\quad U=\frac{\k_1(\gamma-1)\Theta_x}{\gamma R\Theta},
\quad \Theta=\Theta.
\end{equation}
The contact wave $(V, U,\Theta)(x,t)$ solves the compressible Navier-Stokes system \eqref{ns}
time asymptotically, that is,
\begin{equation}\label{contact equation}
\left\{
\begin{array}{ll}
\di V_t-U_x=0,\\
\di U_t+\left(\frac{R\Theta}{V}\right)_x=\mu\left(\frac{U_x}{V}\right)_x+R_1,\\
\di c_{\nu}\Theta_t+p(V,\Theta)U_x=\left(\kappa_1 \frac{\Theta_x}{V}\right)_x+\mu\frac{U^2_x}{V}+R_2,
\end{array}
\right.
\end{equation}
where
\begin{equation}
\di \widetilde R_1=U_t-\mu\left(\frac{U_x}{V}\right)_x,\quad \widetilde R_2=-\mu\frac{U_x^2}{V}.
\label{r1r2}
\end{equation}

For the single viscous contact wave $(V, U,\Theta)$ with arbitrary $\gamma>1$, let
 the perturbation $(\phi,\psi,\zeta)(x,t)$ by
\begin{equation}
(\phi,\psi,\zeta)(x,t)=(v-V,u-U,\t-\Theta)(x,t).
\end{equation}

Then our first main result is as follows: 
\begin{theorem}{\bf (Viscous contact wave)}\label{theorem1}
Suppose that 
$$\alpha=0,\beta>0.$$
For any given left end state $(v_-,u_-,\theta_-)$, suppose that the right end state $(v_+,u_+,\t_+)$ satisfies \eqref{RH}.
Let $(V,U,\Theta)$ be the viscous contact wave defined in \eqref{contact} with strength
$\d=|\t_+-\t_-|$. There exist a function $m(\d)$ satisfying $m(\d)\to +\infty$, as $\d\to 0$ and a small constant $\delta_0$ such that if $\d<\d_0$ and the initial data satisfies
\begin{equation}\label{data}
\left\{
\begin{array}{ll}
\di v_0(x), \t_0(x)\geq m_0^{-1},\quad m_0=:m(\d_0),\\[0.2cm]
\|(v_0(x)-V(x,0),u_0(x)-U(x,0),\t_0(x)-\Theta(x,0))\|_{H^1(\mathbb{R})}\le m_0,
\end{array}
\right.
\end{equation}
then the Cauchy problem \eqref{ns}-\eqref{initial}  admits
a unique global solution $(v,u,\t)$ satisfying
$$
(v-V,u-U,\theta-\Theta)(x,t)\in C\big((0,+\infty);H^1(\mathbb{R})\big);
$$
$$
(v-V)_x(x,t)\in L^2(0,+\infty;L^2\big(\mathbb{R})\big);
$$
$$
(u-U,\t-\Theta)_x(x,t)\in L^2\big(0,+\infty;H^1(\mathbb{R})\big).
$$
Furthermore,
\begin{equation}\label{behavior-1}
\lim_{t\rightarrow+\infty}\sup_{x\in\mathbb{R}}\left|(v-V,u-U,\t-\Theta)(x,t)\right|=0.
\end{equation}

\end{theorem}

\begin{remark}
Theorem 1.1 can be regarded as a natural generalization of Huang-Wang's result(\cite{H-W}) where they considered the constant viscosity case $(\alpha=\beta=0)$ to the degenerate and nonlinear one that $\alpha=0,\beta>0$. It is interesting to study the case that $\alpha>0,\beta>0$, which will be left for future. 
\end{remark}


When the relation \eqref{RH} fails, the basic theory of hyperbolic systems of conservation laws
implies that for any given constant state $(v_-, u_-, \t_-)$ with $v_->0$, $\t_->0$ and $u_-\in\mathbb{R}$,
there exists a suitable neighborhood $\Omega(v_-,u_-,\t_-)$ of $(v_-, u_-, \t_-)$ such that for any
$(v_+, u_+, \t_+)\in\Omega(v_-,u_-,\t_-)$, the Riemann problem of the Euler system \eqref{euler}, \eqref{Riemann}
has a unique solution. In this paper, we only consider the stability of the superposition of the viscous
contact wave and rarefaction waves. In this situation, we assume that 
\begin{equation}\label{zheng}
\di (v_+,u_+,\t_+)\in R_1CR_3(v_-,u_-,\t_-)\subset\Omega(v_-,u_-,\t_-),
\end{equation}
where
\begin{equation*}
\begin{array}{ll}
\di u\geq u_--\int_{v_-}^{e^{\frac{\g-1}{R\g}(s_--s)}v}\lambda_-(\eta,s_-)d\eta,
u\geq u_--\int^v_{e^{\frac{\g-1}{R\g}(s_--s)}v_-}\lambda_+(\eta,s)d\eta\Bigg\},
\end{array}
\end{equation*}
$$ s=\frac{R}{\g-1}\ln\frac{R\t}{A}+R\ln v, $$\\
$$s_{\pm}=\frac{R}{\g-1}\ln\frac{R\t_{\pm}}{A}+R\ln v_{\pm},~~\\
$$\\
$$ \lambda_{\pm}(v,s)=\pm(A\g v^{-\g-1}e^{\frac{\g-1}{R}s})^{1/2},$$
and
$$
 R_1CR_3(v_-,u_-,\t_-)\triangleq\Bigg\{(v,u,\t)\in\Omega(v_-,u_-,\t_-)\Bigg|s\neq s_-.
$$

By the standard argument, there exists a unique pair of points $(v_-^m, u^m, \t_-^m)$
and $(v_+^m, u^m, \t_+^m)$ in $\Omega(v_-,u_-,\t_-)$ satisfying
$$
\frac{R\t_-^m}{v_-^m}=\frac{R\t_+^m}{v_+^m}\triangleq p^m,
$$
such that the points $(v_-^m,u^m,\t_-^m)$ and $(v_+^m,u^m,\t_+^m)$ belong to the 1-rarefaction wave
curve $R_-(v_-,u_-,\t_-)$ and the 3-rarefaction wave curve $R_+(v_+,u_+,\t_+)$, respectively, where
$$
R_{\pm}(v_{\pm},u_{\pm},\t_{\pm})=\left\{(v,u,\t)\in\Omega(v_-,u_-,\t_-)\Bigg|s=s_{\pm},
u=u_{\pm}-\int_{v_{\pm}}^v\lambda_{\pm}(\eta,s_{\pm})d\eta,v>v_{\pm}
\right\}.
$$

We assume $u^m=0$ in what follows without loss of generality.
The 1-rarefaction wave $(v_-^r,u_-^r,\t_-^r)(\frac{x}{t})$ (respectively the 3-rarefaction wave $(v_+^r,u_+^r,\t_+^r)(\frac{x}{t})$)
connecting $(v_-,u_-,\t_-)$ and $(v_-^m,0,\t_-^m)$ (respectively $(v_+^m,0,\t_+^m)$ and $(v_+,u_+,\t_+)$) is the
weak solution of the Riemann problem of the Euler system \eqref{euler} with the following initial Riemann data
\begin{equation}\label{Riemann2}
\di(v_{\pm},u_{\pm},\t_{\pm})(x,0)=\left\{
\begin{array}{ll}
(v_{\pm}^m,0,\t_{\pm}^m), &\quad\pm x<0,\\
(v_{\pm},u_{\pm},\t_{\pm}), &\quad\pm x>0.
\end{array}
\right.
\end{equation}
Since the rarefaction wave $(v_{\pm}^r,u_{\pm}^r,\t_{\pm}^r)$ are not smooth enough solutions, it is
convenient to construct approximate rarefaction wave which is smooth. Motivated by \cite{MN-86}, the smooth solutions of Euler
system \eqref{euler}, $(V_{\pm}^r,U_{\pm}^r,\T_{\pm}^r)$, which approximate $(v_{\pm}^r,u_{\pm}^r,\t_{\pm}^r)$,
are given by
\begin{equation}\label{appro-rare}
\di\left\{
\begin{array}{ll}
\di\lambda_{\pm}(V_{\pm}^r(x,t),s_{\pm})=w_{\pm}(x,t),\\
\di U_{\pm}^r=u_{\pm}-\int_{v_{\pm}}^{V_{\pm}^r(x,t)}\lambda_{\pm}(\eta,s_{\pm})d\eta,\\
\T_{\pm}^r=\t_{\pm}(v_{\pm})^{\g-1}(V_{\pm}^r)^{1-\g},
\end{array}
\right.
\end{equation}
 where $w_-$ (respectively $w_+$) is the solution of the initial problem for the typical Burgers equation:
\begin{equation}\label{burgers}
\di\left\{
\begin{array}{ll}
\di w_t+ww_x=0,\quad(x,t)\in\mathbb{R}\times(0,\infty),\\
\di w(x,0)=\frac{w_r+w_l}{2}+\frac{w_r-w_l}{2}\tanh x,
\end{array}
\right.
\end{equation}
with $w_l=\lambda_-(v_-,s_-)$, $w_r=\lambda_-(v_-^m,s_-)$ (respectively $w_l=\lambda_+(v_+^m,s_+)$, $w_r=\lambda_+(v_+,s_+)$).

Let $(V^{cd},U^{cd},\T^{cd})(x,t)$ be the viscous contact wave constructed in \eqref{1.8} and \eqref{contact}
with $(v_{\pm},u_{\pm},\t_{\pm})$ replaced by $(v_{\pm}^m,0,\t_{\pm}^m)$, respectively.

For subsequent use, we define the strengths of the viscous contact wave and rarefaction waves as follow:
\begin{equation*}
\begin{array}{ll}
\di\d^{r_1}=|v_-^m-v_-|+|0-u_-|+|\t_-^m-\t_-|, \quad\\[3mm] \di\d^{cd}=|\t_+^m-\t_-^m|,\\[3mm]
\di\d^{r_3}=|v_+^m-v_+|+|0-u_+|+|\t_+^m-\t_+|
\end{array}
\end{equation*}
and $\d=\min(\d^{r_1},\d^{cd},\d^{r_3})$. If
\begin{equation}\label{same-order}
\d^{r_1}+\d^{cd}+\d^{r_3}\leq C\d,\quad {\rm as }\quad \d^{r_1}+\d^{cd}+\d^{r_3}\rightarrow 0
\end{equation}
holds for a positive constant $C$, we call the strengths of the wave patterns ``small with the same order". In this case, we have
\begin{equation}
\d^{r_1}+\d^{cd}+\d^{r_3}\leq C|(v_+-v_-,u_+-u_-,\t_+-\t_-,)|.
\end{equation}
Henceforth, we assume \eqref{same-order} holds throughout. We define
\begin{equation}\label{ansatz}
\left(
\begin{array}{l}
\di V\\
\di U \\
\di \T
\end{array}
\right)(x,t) = \left(
\begin{array}{l}
\di V^{cd}+V_-^r+V_+^r\\
\di U^{cd}+U_-^r+U_+^r\\
\di \T^{cd}+\T_-^r+\T_+^r
\end{array}\right)(x,t)-\left(
\begin{array}{l}
\di v_-^m+v_+^m\\
\di \quad~~ 0\\
\di \t_-^m+\t_+^m
\end{array}\right),
\end{equation}
and
$$
(\phi,\psi,\zeta)(x,t)=(v-V,u-U,\t-\T)(x,t).
$$

Then our second main result is as follows:
\begin{theorem}{\bf (Composite waves)}\label{theorem2}Suppose that 
$$\alpha=0,\beta>0.$$
For any given left end state $(v_-,u_-,\theta_-)$, let $(V,U,\Theta)$ be defined in \eqref{ansatz} with strength
satisfying \eqref{same-order}. Then there exists a function $m(\d)$ satisfying $m(\d)\to +\infty$,
as $\d\to 0$ and a small constant $\delta_0$ , such that if $|(v_+-v_-,u_+-u_-,\t_+-\t_-)|<\d_0$ and the initial data satisfies
\begin{equation}\label{data2}
\left\{
\begin{array}{ll}
\di  v_0(x), \t_0(x)\geq m_0^{-1},\quad m_0=:m(\d_0),\\[0.2cm]
\|(v_0(x)-V(x,0),u_0(x)-U(x,0),\t_0(x)-\Theta(x,0))\|_{H^1(\mathbb{R})}\le m_0,
\end{array}
\right.
\end{equation}
then the Cauchy problem \eqref{ns}-\eqref{initial} admits
a unique global solution $(v,u,\t)$ satisfying
$$
(v-V,u-U,\theta-\Theta)(x,t)\in C\big((0,+\infty);H^1(\mathbb{R})\big);
$$
$$
(v-V)_x(x,t)\in L^2(0,+\infty;L^2\big(\mathbb{R})\big);
$$
$$
(u-U,\t-\Theta)_x(x,t)\in L^2\big(0,+\infty;H^1(\mathbb{R})\big),
$$
and
\begin{equation}\label{wt-2}
\lim_{t\rightarrow+\infty}\sup_{x\in\mathbb{R}}\left|(v-V,u-U,\t-\Theta)(x,t)\right|=0,
\end{equation}
where the $(v_-^r,u_-^r,\t_-^r)(x,t)$ and $(v_+^r,u_+^r,\t_+^r)(x,t)$ are the 1-rarefaction and
3-rarefaction waves uniquely determined by \eqref{euler}, \eqref{Riemann2}, respectively.
\end{theorem}

\begin{remark}
By $(iv)$ of Lemma \ref{rare-pro}, Theorem \ref{theorem2} implies
\begin{equation*}\label{behavior-2}
\lim_{t\rightarrow+\infty}\sup_{x\in\mathbb{R}}
\left(
\begin{array}{l}
|(v-v_-^r-V^{cd}-v_+^r+v_-^m+v_+^m)(x,t)|\\[2mm]
\quad\quad|(u-U^{cd}-u_-^r-u_+^r)(x,t)|\\[2mm]
|(\t-\t_-^r-\T^{cd}-\t_+^r+\t_-^m+\t_+^m)(x,t)|\\[2mm]
\end{array}
\right)=0.
\end{equation*}
\end{remark}

We now explain the main analysis of this paper. Compared with  the constant viscosity case$(\alpha=\beta=0)$(\cite{H-W}), the main difficulty comes from the degeneracy and nonlinearity of the heat conductivity due to the fact that $\beta>0$. To overcome this difficulty, we utilize the structure of wave patterns to essentially control the terms involving the derivative of perturbation around the wave patterns.  Firstly, we make some slight modifications to the ideas in \cite{Li-Xu,jiang} to get the time-independent lower and upper bounds of $v$,  provided the strengths of the waves are suitably small. Next, combining  the uniform bound of $v$, the basic energy estimate(see (\ref{basic})) and cut-off techniques, we are just motivated by \cite{Li-Xu} to consider the cases $\beta\in (0,1)$ and $\beta \in (1,+\infty) $ separately. This approach yields the key issue of $\psi_x$ and $\t^{1/2}\zeta_x$ on the $L^2(\mathbb{R}\times (0,T))$-norm (see (\ref{new-nergy}) and (\ref{beta-new})). Then, we modify some ideas of \cite{Li-Xu,H-W} to obtain the estimates on the $L^2(\mathbb{R}\times (0,T))$-norm of $\zeta_x,\phi_x$ and $\psi_{xx}^2$ with suitably small $\delta_0$(see lemma \ref{high-order}). Finally, to resolve the difficulty arising from cross term, we observe the simple but crucial fact (see (\ref{k11})) 
$$\int\left(\frac{\t^\beta\zeta_x}{v}\right)_x\t^\beta\zeta_tdx =\int\left(\frac{\t^\beta\zeta_x}{v}\right)\left((\t^\beta\zeta_x)_t + \beta\theta^{\beta-1}(\zeta_t\T_x-\zeta_x\T_t)\right)dx, $$
which combined with (\ref{int-t+}
)-(\ref{k5}), (\ref{basic}), (\ref{theta^beta}) and the Gronwall's inequality leads to the $L^\infty(0,T;L^2(\mathbb{R}))$-norm of $\t^\beta\zeta_x$. With these weighted estimates on the perturbation around the wave patterns at hand, we derive the upper and lower bounds of $\t$. These are the key to the proof, and the whole procedure will be carried out in the next section.


\

This paper is organized as follows. In the next section, we collect some useful lemmas and
fundamental facts concerning the viscous contact wave as well as rarefaction waves.
The main proof of Theorem \ref{theorem1} and \ref{theorem2} are completed in Section 4 and Section 5, respectively.

\

%
%
%


\

\section{Preliminaries}
\setcounter{equation}{0}

In this section, we first consider a lemma for the properties of the viscous contact wave $(V,U,\Theta)$ defined by \eqref{contact} are useful  in the following sections.
\begin{lemma}(lemma2.1 in \cite{H-W})\label{decay}
Assume that $\delta=|\t_+-\t_-|\leq \d_0$ for a small positive constant $\d_0$. Then the viscous contact wave $(V,U,\Theta)$ has the following properties:
$$
|V-v_{\pm}|+|\Theta-\t_{\pm}|\leq O(1)\d e^{-\frac{c_1x^2}{1+t}},
$$
$$
|\partial^k_x V|+|\partial^{k-1}_x U|+|\partial^k_x\Theta|\leq
O(1)\delta(1+t)^{-\frac{k}{2}}e^{-\frac{c_1x^2}{1+t}}, \quad k\geq 1.
$$

\end{lemma}
Combining (\ref{r1r2}) we have
\begin{equation}
\di \widetilde R_1=O(1)\d(1+t)^{-\frac{3}{2}}e^{-\frac{c_1x^2}{1+t}},
\quad \widetilde R_2=O(1)\d(1+t)^{-2}e^{-\frac{c_1x^2}{1+t}}.
\end{equation}

The next  two lemmas  play important roles to obtain the basic energy
estimate.
\begin{lemma}(Lemma 1 in \cite{Huang-Li-matsumura})\label{hlm}
For $0<T\leq+\infty$, suppose that $h(x,t)$ satisfies
$$
h\in L^{\infty}(0,T;L^2(\mathbb{R})),\quad h_x\in L^{2}(0,T;L^2(\mathbb{R})),\quad h_t\in L^{2}(0,T;H^{-1}(\mathbb{R})).
$$
Then the following estimate holds 
\begin{equation}
\int_0^T\int h^2w^2dxdt\leq 4\pi\|h(0)\|^2+4\pi\alpha^{-1}\int_0^T\|h_x\|^2dt
+8\alpha\int_0^T<h_t,hg^2>_{H^{-1}\times H^1}dt
\end{equation}
for $\alpha>0$, and
\begin{equation}
	w(x,t)=(1+t)^{-\frac{1}{2}}\exp\left(-\frac{\alpha x^2}{1+t}\right),\quad
	g(x,t)=\int_{-\infty}^x w(y,t)dy,
	\label{w}
\end{equation}

\end{lemma}
and $<\cdot,\cdot>$ denotes the dual product between $H^{-1}$ and $H^1$.\\

\begin{lemma}\label{lemma5}(lemma 5 in\cite{Huang-Li-matsumura} )
For $\alpha\in(0,\frac{c_1}{4}]$ and $w$ defined by (\ref{w}) , there exists some positive
constant $C$ depending on $\alpha$ such that the following estimate holds
\begin{equation}\label{important}
\int_0^t\int(\phi^2+\psi^2+\zeta^2)w^2 dxds\leq C\left(
1+\int_0^t\int (\phi_x^2+\psi_x^2+\zeta_x^2)dxds
\right). 
\end{equation}
\end{lemma}

Next, we brevity state the following properties of the solution to the problem \eqref{burgers}, the proofs can be found in \cite{MN-86}.

\begin{lemma}\label{bur-pro}
For given $w_l\in\mathbb{R}$ and $\bar w>0$, let $w_r\in\{0<\tilde w\triangleq w-w_l<\bar w\}$.
Then the problem \eqref{burgers} has a unique smooth global solution in time satisfying the following properties.
\begin{itemize}
\item[(1)] $w_l<w(x,t)<w_r$, $w_x>0$ $(x\in\mathbb{R},t>0)$.
\item[(2)] For $p\in[1,\infty]$, there exists some positive constant $C=C(p,w_l,\bar w)$ such that
for $\tilde w\geqq0$ and $t\geqq0$,
$$
\|w_x(t)\|_{L^p}\leq C\min\{\tilde w,\tilde w^{1/p}t^{-1+1/p}\},\quad
\|w_{xx}(t)\|_{L^p}\leq C\min\{\tilde w,t^{-1}\}.
$$
\item[(3)] If $w_l>0$, for any $(x,t)\in(-\infty,0]\times[0,\infty)$,
$$
|w(x,t)-w_l|\leq\tilde we^{-2(|x|+w_lt)},\quad
|w_x(x,t)|\leq2\tilde we^{-2(|x|+w_lt)}.
$$
\item[(4)] If $w_r<0$, for any $(x,t)\in[0,\infty)\times[0,\infty)$,
$$
|w(x,t)-w_r|\leq\tilde we^{-2(x+|w_r|t)},\quad
|w_x(x,t)|\leq2\tilde we^{-2(x+|w_r|t)}.
$$
\item[(5)] For the Riemann solution $w^r(x/t)$ of the scalar equation \eqref{burgers} with the Riemann
initial data
\begin{equation*}
w(x,0)=\left\{
\begin{array}{ll}
w_l,&\quad x<0,\\
w_r,&\quad x>0,
\end{array}
\right.
\end{equation*}
we have
$$
\lim_{t\rightarrow+\infty}\sup_{x\in\mathbb{R}}|w(x,t)-w^r(x/t)|=0.
$$
\end{itemize}

\end{lemma}

We finally divide $\mathbb{R}\times(0,t)$ into
three parts, that is $\mathbb{R}\times(0,t)=\Omega_-\cup\Omega_c\cup\Omega_+$ with
$$
\Omega_{-}=\big\{(x,t)\big|2x<\lambda_{-}(v_{-}^m,s_{-})t\big\},
$$
$$
\Omega_{c}=\big\{(x,t)\big|\lambda_-(v_-^m,s_-)t\leq2x\leq\lambda_{+}(v_{+}^m,s_{+})t\big\},
$$
$$
\Omega_{+}=\big\{(x,t)\big|2x>\lambda_{+}(v_{+}^m,s_{+})t\big\}.
$$
We have
\begin{lemma}(lemma2.5 in \cite{H-W})\label{rare-pro}
For any given left end state $(v_-,u_-,\t_-)$, suppose that \eqref{zheng}
and \eqref{same-order} hold. Then the smooth rarefaction waves $(V_{\pm}^r,U_{\pm}^r,\T_{\pm}^r)$
constructed in \eqref{appro-rare} and the viscous contact wave $(V^{cd},U^{cd},\T^{cd})$ constructed in \eqref{contact} satisfying
the following:
\begin{itemize}
\item[(1)] $(U_{\pm}^r)_x\geq0$,~$(x\in\mathbb{R},t>0)$.
\item[(2)] For $p\in[1,\infty]$, there exists a positive constant $C=C(v_-,u_-,\t_-,\d)$
such that for $\d$ satisfying \eqref{same-order},
$$
\|\big((V_{\pm}^r)_x,(U_{\pm}^r)_x,(\T_{\pm}^r)_x\big)(t)\|_{L^p}
\leq C\min\Big\{\d,\d^{1/p}t^{-1+1/p}\Big\}
$$
and
$$
\|\big((V_{\pm}^r)_{xx},(U_{\pm}^r)_{xx},(\T_{\pm}^r)_{xx}\big)(t)\|_{L^p}
\leq C\min\Big\{\d,t^{-1}\Big\}.
$$
\item[(3)] There exists a positive constant $C=C(v_-,u_-,\t_-,\d)$
such that for
$$
c_0=\frac{1}{10}\min\Big\{|\lambda_-(v_-^m,s_-)|,\lambda_+(v_+^m,s_+),c_1\lambda_-^2(v_-^m,s_-),
c_1\lambda_+^2(v_+^m,s_+),1\Big\},
$$
we have 
\begin{equation*}
	(U_{\pm}^r)_x+|(V_{\pm}^r)_x|+|V_{\pm}^r-v_{\pm}^m|+|(\T_{\pm}^r)_x|+|\T_{\pm}^r-\t_{\pm}^m|
	\leq C\d e^{-c_0(|x|+t)}, \quad in \Omega_c,
\end{equation*}
\begin{equation*}
	\left\{
	\begin{array}{ll}
		|V^{cd}-v_{\mp}^m|+|V^{cd}_x|+|\T^{cd}-\t_{\mp}^m|+|U^{cd}_x|+|\T^{cd}_x|
		\leq C\d e^{-c_0(|x|+t)},\\
		(U_{\pm}^r)_x+|(V_{\pm}^r)_x|+|V_{\pm}^r-v_{\pm}^m|+|(\T_{\pm}^r)_x|+|\T_{\pm}^r-\t_{\pm}^m|
		\leq C\d e^{-c_0(|x|+t)},
	\end{array}
	\quad in \Omega_{\mp}.
	\right.
\end{equation*}
\item[(4)] For the rarefaction waves $(v_{\pm}^r,u_{\pm}^r,\t_{\pm}^r)(x/t)$ determined by
\eqref{euler} and \eqref{Riemann2}, it holds that 
$$
\lim_{t\rightarrow+\infty}\sup_{x\in\mathbb{R}}
\big|(V_{\pm}^r,U_{\pm}^r,\T_{\pm}^r,)(x,t)-(v_{\pm}^r,u_{\pm}^r,\t_{\pm}^r,)(x/t)\big|=0.
$$
\end{itemize}

\end{lemma}


\section{Proof of Theorem \ref{theorem1}}
\setcounter{equation}{0}
For equations \eqref{contact equation}, we consider Cauchy problem \eqref{ns} and  \eqref{initial} yields
\begin{equation}\label{perturb}
\left\{
\begin{array}{ll}
\di \phi_t-\psi_x=0,\\
\di \psi_t+\left(p-p_+\right)_x=\tilde{\mu}\left(\frac{u_x}{v}-\frac{U_x}{V}\right)_x-\widetilde R_1,\\
\di c_{\nu}\zeta_t+pu_x-p_+U_x=\tilde{\kappa}\left(\frac{\t^{\beta}\t_x}{v}-\frac{\Theta^{\beta}\Theta_x}{V}\right)_x
+\tilde{\mu}\left(\frac{u_x^2}{v}-\frac{U^2_x}{V}\right)-\widetilde R_2,\\
\di (\phi,\psi,\zeta)(x,0)=(\phi_0,\psi_0,\zeta_0)(x),\quad x\in \mathbb{R}.
\end{array}
\right.
\end{equation}
For the solution $(\phi,\psi,\zeta)$, we define the solution space $X([0,+\infty))$,
\begin{equation*}
\begin{array}{l}
\di X([0,T])=\Big\{(\phi,\psi,\zeta)\Big| v, ~\theta \geq M^{-1},
 \quad \sup_{0\leq t\leq T}\|(\phi,\psi,\zeta)\|_{H^1}\leq M
 \Big\}
\end{array}
\end{equation*}
for some $0<T\leq +\infty$, where the constants $M$ will be determined later.
Since the local existence of the solution is well known (see \cite{huang-ma-shi}), we can prove the global existence part of  Theorem \ref{theorem1} by the following a priori estimates.

\begin{proposition}\label{prop}
\textbf{(A priori estimates)} Assume that the conditions of Theorem \ref{theorem1} hold, then there exists
a positive constant $\d_0$ such that if $\d< \d_0$ and $(\phi,\psi,\zeta)\in X([0,T])$,
\begin{equation}
\begin{array}{l}\label{p-1}
\di \sup_{0\leq t\leq T}\|(\phi,\psi,\zeta)(t)\|_{H^1}^2+\int_0^T(\|\phi_x\|^2+\|(\psi_x,\zeta_x)\|^2_{H^1})ds
\leq C_0,
\end{array}
\end{equation}
where $C_0$ denotes a constant depending only on $\tilde{\mu}$, $\tilde{\kappa}$, $R$, $c_{\nu}$, $v_{\pm}$, $u_{\pm}$, $\t_{\pm}$ and  $m_0$.

\end{proposition}

Once Proposition \ref{prop} is proved, we can extend the unique local solution $(u,v,\t)$ which can be
obtained as in \cite{huang-ma-shi} to $T=\infty$. Estimate \eqref{p-1}
and the equations \eqref{perturb} (respectively \eqref{perturb2}) imply that
\begin{equation}
\di\int_0^{+\infty}\left(\|(\phi_x, \psi_x,\zeta_x)(t)\|^2+\left|\frac{d}{dt}\|(\phi_x, \psi_x,\zeta_x)(t)\|^2\right|\right)dt< \infty,
\end{equation}
which together with \eqref{p-1} and the Sobolev's inequality, easily leads to the large time behavior of the solutions  \eqref{behavior-1} (respectively \eqref{wt-2}).

 Proposition \ref{prop} will be finished by the following lemmas. 

\begin{lemma}\label{basic-lemma}
There exist some positive constant $C_0$  and $\d_0$ such that if $\d<\d_0$, it holds that
\begin{equation}\label{basic}
\begin{array}{ll}
\di \Big\|\left(\psi,\sqrt{\Phi\left(\frac{v}{V}\right)},\sqrt{\Phi\left(\frac{\t}{\Theta}\right)}\right)(t)\Big\|^2
+\int_0^t\int\left(\frac{\psi_x^2}{\t v}+\frac{\t^\beta\zeta_x^2}{\t^2 v}\right)dxds
 \leq C_0,
\end{array}
\end{equation}
where define 
\begin{equation}
	V(t) = \int\left(\frac{\psi_x^2}{\t v}+\frac{\t^\beta\zeta_x^2}{\t^2 v}\right)dx
	\label{v(t)}.
\end{equation}

\end{lemma}

\textbf{Proof}: The proof of the Lemma \ref{basic-lemma} consists of two steps.

\underline{ Step 1. }\quad Similar to \cite{huang-ma-shi}, multiplying $\eqref{ns}_1$ by
$-R\Theta(v^{-1}-V^{-1})$, $\eqref{ns}_2$ by $\psi$ and $\eqref{ns}_3$ by $\zeta\theta^{-1}$, then
adding the resulting equations together yield
\begin{equation}
\begin{array}{ll}
\di\left(\frac{\psi^2}{2}+R\Theta\Phi\left(\frac{v}{V}\right)+c_{\nu}\Theta\Phi\left(\frac{\theta}{\Theta}\right)\right)_t
+\frac{\tilde{\mu}\Theta}{\theta v}\psi_x^2+\frac{\tilde{\kappa}	\t^\beta\Theta}{\theta^2 v}\zeta_x^2
+H_{x}+Q=-\psi \widetilde R_1-\frac{\zeta}{\theta}\widetilde R_2,
\end{array}
\label{ba1}
\end{equation}
where
$$
\Phi(z)=z-\ln z-1, \quad z>0,
$$
\begin{equation}
	\label{H}
\di H=(p-p_+)\psi-\tilde{\mu}\left(\frac{u_x}{v}-\frac{U_x}{V}\right)\psi-\frac{\tilde{\k}\zeta}{\theta}\left(\frac{\t^\beta\t_x}{v}-\frac{\T^\beta\Theta_x}{V}\right),
\end{equation}
\begin{equation*}
\begin{array}{ll}
\di Q&\di=p_+\Phi\left(\frac{V}{v}\right)U_x+\frac{p_+}{\gamma-1}\Phi\left(\frac{\Theta}{\theta}\right)U_x
+\tilde{\mu}\left(\frac{1}{v}-\frac{1}{V}\right)\psi_xU_x\\[3mm]
&\di -\frac{\zeta}{\theta}(p_+-p)U_x-\frac{\tilde{\k}\t^\beta\Theta_x}{\theta^2v}\zeta\zeta_x
-\frac{\tilde{\k}\Theta\Theta_x}{\theta^2vV}(\Theta^\beta v-\theta^\beta V)\zeta_x\\[3mm]
&\di+\frac{\tilde{\k}\Theta_x^2}{\t^2vV}(\Theta^\beta v-\theta^\beta V)\zeta
-\frac{2\tilde{\mu} U_x}{\theta v}\psi_x\zeta+\frac{\tilde{\mu} U_x^2}{\theta vV}\zeta\phi.
\end{array}
\end{equation*}
Since
\begin{equation}
	\label{Q}
	\begin{array}{ll}
\di Q&\di\leq \frac{\tilde{\mu}\Theta}{4\theta v}\psi_x^2+\frac{\tilde{\kappa}\t^\beta\Theta}{8\theta^2v}\zeta_x^2+\frac{\tilde{\kappa}\Theta\Theta_x}{\t^2vV}\zeta_x \left(\Theta^\beta\zeta + V(\Theta^\beta-\t^\beta)\right)   \\[3mm]
&\di +\frac{\tilde{\kappa}\Theta^2_x}{\t^2vV}\zeta\left(\Theta^\beta\phi+V(\Theta^\beta-\theta^\beta)\right)+ C(M)(\psi^2+\zeta^2)(|U_x|+\Theta_x^2)\\[3mm]
&\di\leq \frac{\tilde{\mu}\Theta}{4\theta v}\psi_x^2+\frac{\tilde{\kappa}\t^\beta\Theta}{8\theta^2v}\zeta_x^2+\frac{\tilde{\kappa}\Theta\Theta_x}{\t^2vV}\zeta_x \left(\Theta^\beta\zeta + C(M)|\zeta|\right)  \\[3mm]
&\di +\frac{\tilde{\kappa}\Theta^2_x}{\t^2vV}\zeta\left(\Theta^\beta\phi+C(M)|\zeta|\right)+ C(M)(\psi^2+\zeta^2)(|U_x|+\Theta_x^2)\\[3mm]
&\di\leq \frac{\tilde{\mu}\Theta}{4\theta v}\psi_x^2+\frac{\tilde{\kappa}\t^\beta\Theta}{4\theta^2v}\zeta_x^2
+C(M)(\phi^2+\zeta^2)(|U_x|+\Theta_x^2),
	\end{array}
\end{equation}
where $C(M)$ denotes a constant depending on $M$. Recalling Lemma \ref{decay}, we have
\begin{equation}\label{R1phi}
\begin{array}{ll}
\di \left|\int_0^t\int\widetilde R_1\psi dxds\right|&\di\leq O(1)\d \int_0^t\int (1+s)^{-3/2}e^{-\frac{c_1x^2}{1+s}}|\psi|dxds\\[3mm]
&\di \leq O(1)\d\int_0^t(1+s)^{-5/4}\|\psi\|ds\\[3mm]
&\di \leq O(1)\d\int_0^t(1+s)^{-5/4}\|\psi\|^2ds+O(1)\d
\end{array}
\end{equation}
and
\begin{equation}\label{R2zeta}
\begin{array}{ll}
\di \left|\int_0^t\int\widetilde R_2\frac{\zeta}{\t} dxds\right|
\leq C(M)\d\int_0^t(1+s)^{-7/4}\left\|\sqrt{\Phi\left(\frac{\t}{\T}\right)}\right\|^2ds+C(M)\d.
\end{array}
\end{equation}
Then, integrating \eqref{ba1} over $\mathbb{R}\times(0,t)$, choosing $\a=\frac{c_1}{4}$ in Lemma \ref{hlm}
and $\d$ suitable small, it follows from Lemma \ref{hlm}-\ref{lemma5}, (\ref{Q})-(\ref{R2zeta}) and Gronwall's inequality that
\begin{equation}\label{ba2}
\begin{array}{ll}
\di \Big\|\left(\psi,\sqrt{\Phi\left(\frac{v}{V}\right)},\sqrt{\Phi\left(\frac{\t}{\Theta}\right)}\right)(t)\Big\|^2
+\int_0^t\int\left(\frac{\psi_x^2}{\t v}+\frac{\t^\beta\zeta_x^2}{\t^2 v}\right)dxds\\[3mm]
\di
\leq C_0+ C(M)\d \int_0^t\int(1+s)^{-1}(\phi^2+\zeta^2)e^{-\frac{c_1x^2}{1+s}}dxds\\
\di \leq C_0+ C(M)\d \int_0^t\int\frac{\t\phi_x^2}{v^3}dxds.
\end{array}
\end{equation}

\underline{ Step 2. }\quad  Following the approach of \cite{MN-92}, we define a new
variable $\di\tilde v=\frac{v}{V}$. Then $\eqref{perturb}_2$ can be rewritten by the
new variable as
\begin{equation}\label{newperturb}
\di \left(\tilde{\mu}\frac{\tilde v_x}{\tilde v}-\psi\right)_t-p_x=\widetilde R_1.
\end{equation}
Multiplying \eqref{newperturb} by $\di\frac{\tilde v_x}{\tilde v}$, we obtain
\begin{equation}\label{phi-x1}
\begin{array}{ll}
\di \left(\frac{\tilde{\mu}}{2}\left(\frac{\tilde v_x}{\tilde v}\right)^2-\psi\frac{\tilde v_x}{\tilde v}\right)_t
+\left(\psi\frac{\tilde v_t}{\tilde v}\right)_x+\frac{R\t}{v}\left(\frac{\tilde v_x}{\tilde v}\right)^2
-\frac{R}{v}\zeta_x\frac{\tilde v_x}{\tilde v}\\[3mm]
\di\quad +\frac{R\t}{v}\left(\frac{1}{\Theta}-\frac{1}{\t}\right)\Theta_x\frac{\tilde v_x}{\tilde v}
=\frac{\psi_x^2}{v}+\psi_xU_x\left(\frac{1}{v}-\frac{1}{V}\right)+\widetilde R_1\frac{\tilde v_x}{\tilde v}.
\end{array}
\end{equation}
It follows from Cauchy's inequality yields that
\begin{equation}
	\label{4.12}
\begin{array}{ll}
\di \left|\frac{R}{v}\zeta_x\frac{\tilde v_x}{\tilde v}\right|
+\left|\psi_xU_x\left(\frac{1}{v}-\frac{1}{V}\right)\right|+\frac{\psi_x^2}{v}\\[3mm]
\di \leq \frac{R\t}{4v}\left(\frac{\tilde v_x}{\tilde v}\right)^2+
C(M)\left(\frac{\t^\beta\zeta_x^2}{\t^2 v}+\frac{\psi_x^2}{\t v}\right)
+C(M)\phi^2U_x^2,
\end{array}
\end{equation}
and
\begin{equation}
\begin{array}{ll}
	\label{4.13}
\di \left|\frac{R\t}{v}\left(\frac{1}{\Theta}-\frac{1}{\t}\right)\Theta_x\frac{\tilde v_x}{\tilde v}\right|
+\left|\widetilde R_1\frac{\tilde v_x}{\tilde v}\right|\\[3mm]
\di\leq \frac{R\t}{4v}\left(\frac{\tilde v_x}{\tilde v}\right)^2
+C(M)(\zeta^2\Theta_x^2+\widetilde R_1^2).
\end{array}
\end{equation}
Note that
\begin{equation}
\begin{array}{ll}
	\label{4.14}
\di \frac{\phi_x^2}{2v^2}-C(M)\phi^2\Theta_x^2\leq \left(\frac{\tilde v_x}{\tilde v}\right)^2
\leq \frac{\phi_x^2}{v^2}+C(M)\phi^2\Theta_x^2.
\end{array}
\end{equation}
Integrating \eqref{phi-x1} over $\mathbb{R}\times(0,t)$, we get after using (\ref{4.12})-(\ref{4.14}) that
\begin{equation}
\begin{array}{ll}
\di\int\frac{\phi_x^2}{v^2}dx+\int_0^t\int\frac{\t\phi_x^2}{v^3} dxds\leq C_0+C\|\psi\|^2\\[3mm]
\di\quad+C(M)\d^2\left\|\sqrt{\Phi\left(\frac{v}{V}\right)}\right\|^2
+C(M)\int_0^t\int\left(\frac{\psi_x^2}{\t v}+\frac{\t^\beta\zeta_x^2}{\t^2 v}\right) dxds\\[3mm]
\di\quad +C(M)\int_0^t\int(\phi^2+\zeta^2)(U_x^2+\Theta_x^2) dxds.
\end{array}
\end{equation}
By Lemma \ref{lemma5}, \eqref{ba2} and choosing $\d$ suitable small, we have
\begin{equation}\label{phi-x3}
\begin{array}{ll}
\di \int\frac{\phi_x^2}{v^2}dx+\int_0^t\int\frac{\theta\phi_x^2}{v^3}dxds
 \leq C_0+C(M)\int_0^t\int\left(\frac{\psi_x^2}{\t v}+\frac{\t^\beta\zeta_x^2}{\theta^2v}\right) dxds.
\end{array}
\end{equation}
Putting (\ref{phi-x3}) into (\ref{ba2}) and choosing $\delta$ suitable small, the proof of Lemma \ref{basic-lemma} is completed.

\hfill $\Box$

Lemma \ref{basic-lemma} utilizes the smallness of $\delta$ to ensure that the basic energy (\ref{basic}) only depends on the initial data. Building upon the basic energy estimate, we shall establish uniform upper and lower bounds for both the specific volume $v$ and the absolute temperature $\theta$.This determines the requisite smallness of the parameter $\delta$. That is why need the initial condition (\ref{data}). 

To obtain the uniform lower and upper bounds of the specific volume $v(x,t)$, we present the following two lemmas, whose proofs can be found in \cite{H-W}(lemma 3.2) and \cite{jiang}(lemma2.3)
\begin{lemma}{\label{huang-lemma}}
Let $\a_1$, $\a_2$ be the two positive roots of the equation $y-\ln y-1=C_0$ and the constant $C_0$ be the same in \eqref{basic}. Then
\begin{equation}\label{root}
\di \a_1\leq \int_k^{k+1}\tilde v(x,t)dx, \quad\int_k^{k+1}\tilde \t(x,t)dx\leq \a_2,\quad t\geq 0,
\end{equation}
and for each $t\geq 0$ there are points $a_k(t)$, $b_k(t)\in[k,k+1]$ such that
\begin{equation}\label{root2}
\alpha_1\leq \tilde v(a_k(t),t), \tilde \t(b_k(t),t)\leq \alpha_2,\quad t\geq 0,
\end{equation}
where $\tilde v=\frac{v}{V}$, $\tilde \t=\frac{\t}{\Theta}$, and $k=0,\pm1,\pm2,\cdots$.
\end{lemma}


\hfill $\Box$

\begin{lemma}\label{jiang-lemma}
For each $x\in[k,k+1]$, $k=0,\pm1,\pm2,\cdots$, it follows from \eqref{ns}$_2$ that
\begin{equation}\label{v-repre}
\di v(x,t)=B(x,t)Y(t)+\frac{R}{\tilde{\mu}}\int_0^t\frac{B(x,t)Y(t)}{B(x,s)Y(s)}\t(x,s) ds,
\end{equation}
where
\begin{equation}\label{B}
\di B(x,t)=v_0(x)\exp\left(\frac{1}{\tilde{\mu}}\int_x^{+\infty}\big(u_0(y)-u(y,t)\big)D(y)dy\right),
\end{equation}
\begin{equation}\label{Y}
\di Y(t)=\exp\left(\frac{1}{\tilde{\mu}}\int_0^t\int_{k+1}^{k+2}\s(y,s)dyds\right),
\end{equation}
\begin{equation}\label{sigma}
\di \sigma(x,t)=\left(\tilde{\mu}\frac{u_x}{v}-R\frac{\t}{v}\right)(x,t),
\end{equation}
and
\begin{equation}\label{beta}
\di D(x)=\left\{
\begin{array}{ll}
1,& x\leq k+1,\\
k+2-x, & k+1\leq x\leq k+2,\\
0, & x\geq k+2.
\end{array}
\right.
\end{equation}

\end{lemma}
By Cauchy's inequality and \eqref{basic}, we have
\begin{equation}\label{b-bound}
\di \underline{B}(C_0)\leq B(x,t)\leq \overline{B}(C_0), \quad \forall x\in[k,k+1],\quad t\geq 0,
\end{equation}
where $\underline{B}(C_0)$, $\overline{B}(C_0)$ are two constants depending on $C_0$.

With lemma \ref{huang-lemma} and lemma \ref{jiang-lemma} at hand, we are in a position to prove the time-independent upper and lower bounds of $v$.

\begin{lemma}\label{v-bound}
There are two positive constants $\underline{v}(C_0)$, $\bar{v}(C_0)$  such that
\begin{equation}\label{v}
\di \underline{v}(C_0) \leq v(x,t)\leq\bar{v}(C_0),\quad \forall x\in\mathbb{R},\quad t\geq0,
\end{equation}
where $\underline{v}(C_0)$, $\bar{v}(C_0)$  depending on $C_0$, independent of $x$, $t$.

\end{lemma}

\textbf{Proof}:From now on, we always assume that $\t_-<\t_+$ for convenient. So from the properties of
viscous contact wave, we have $\t_-<\T(x,t)<\t_+$ and $v_-<V(x,t)<v_+$.

 It follows from Jensen's inequality,  (\ref{root}), (\ref{root2}) and choose $\delta$ suitable small, we see that for $t\geq s \geq 0$,
\begin{equation}
	\begin{array}{ll}
		\label{int-t/v}
	\di	&\di-\int_{s}^{t}\int_{k+1}^{k+2}\frac{\t}{v}dxd\tau
	\di \leq -C\int_{s}^{t}\int_{k+1}^{k+2}\frac{\tilde{\t}}{\tilde{v}}dxd\tau\\[3mm] 
	&\di\quad\leq C\int_{s}^{t}\int_{k+1}^{k+2}\left(\left(-\tilde{\t}+\tilde \t(b_k(t),t)\right)_+-\tilde \t(b_k(t),t)\right)\frac{1}{\tilde{v}}dxd\tau\\[3mm]
	&\di\quad\leq C\int_{s}^{t}\left(\max_{x\in[k+1,k+2]}(\tilde \t^{\beta/2}(b_k(t),t)-\tilde{\t}^{\beta/2})_+-\tilde \t(b_k(t),t)\right)\left(\int_{k+1}^{k+2}\tilde{v}dx\right)^{-1}d\tau\\[3mm]
	&\di\quad\leq C\int_{s}^{t}\left(\int_{k+1}^{k+2}\tilde{\t}^{\beta/2-1}|\tilde{\t}_x|dx - 2C^{-1}\right)d\tau\\[3mm]
	&\di\quad\leq C\int_{s}^{t}\left(\int_{k+1}^{k+2} \left(\frac{\t^{\beta/2-1}\zeta_x}{\T^{\beta/2}}-\frac{\t^{\beta/2-1}\T_x\zeta}{\Theta^{\beta/2+1}}\right)dx-2C^{-1} \right)d\tau \\[3mm]
	&\di\quad\leq C\int_{s}^{t}\left(\left(\int_{k+1}^{k+2}\frac{\t^\beta\zeta_x^2}{v\t^2}+\frac{\t^\beta\T_x^2\zeta^2}{\t^2\T^2v}dx\right)^{1/2}\left(\int_{k+1}^{k+2}\frac{v}{\T^\beta}dx\right)^{1/2}-2C^{-1}\right)d\tau\\[3mm]
	&\di\quad\leq C\int_{s}^{t}\left(C(\varepsilon)\int_{k+1}^{k+2}\frac{\t^\beta\zeta_x^2}{v\t^2}+\frac{\t^\beta\T_x^2\zeta^2}{\t^2\T^2v}dx + C\varepsilon(v_+\alpha_2)-2C^{-1}\right)d\tau\\[3mm]
	&\di\quad\leq C+ C(M)\int_{s}^{t}\int_{k+1}^{k+2}\zeta^2\T_x^2dxd\tau-\frac{t-s}{C_0}\\[3mm]
	&\di\quad\leq C_0-\frac{t-s}{C_0}.
	\end{array}
\end{equation}

Integrating (\ref{sigma}) over $(k+1,k+2)\times(s,t)$ gives
\begin{equation}\label{int-sigma}
\begin{array}{ll}
\di&\di\int_s^t\int_{k+1}^{k+2}\sigma(x,\tau)dxd\tau=
\int_s^t\int_{k+1}^{k+2}\left(\tilde{\mu}\frac{\psi_x}{v}-R\frac{\t}{v}+\tilde{\mu}\frac{U_x}{v}\right)(x,\tau)dxd\tau\\[3mm]
&\di\quad\leq C\int_s^t\int_{k+1}^{k+2}\frac{\psi_x^2}{\t v}dxd\tau-\frac{R}{2}\int_s^t\int_{k+1}^{k+2}\frac{\t}{v}dxd\tau
+\tilde{\mu}\int_s^t\int_{k+1}^{k+2}\frac{U_x}{v}dxd\tau\\[3mm]
&\di\quad\leq C-\frac{R}{2}\int_s^t\int_{k+1}^{k+2}\frac{\t}{v}dxd\tau + C(M)\int_s^t\int_{k+1}^{k+2} |U_x| dxd\tau\\[3mm]
&\di\quad\leq C+\frac{R}{2}(C_0-\frac{t-s}{C_0})+C(M)\d\\[3mm]
&\di\quad\leq C_0-\frac{t-s}{C_0},\\[3mm]
\end{array}
\end{equation}
where in the third inequality we have used (\ref{int-t/v}). Combining (\ref{int-sigma}) and the definition of $Y(t)$ gives
\begin{equation}\label{Y-e}
\di 0\leq Y(t)\leq C_0e^{-t/C_0},\quad
\frac{Y(t)}{Y(s)}\leq C_0e^{-(t-s)/C_0},
\end{equation}
which together with \eqref{v-repre} and \eqref{b-bound} gives
\begin{equation}\label{v-repre2}
\di v(x,t)\leq C_0e^{-(t-s)/C_0}+C_0\int_0^t\t(x,s)e^{-(t-s)/C_0} ds.
\end{equation}
Then, it follows from (\ref{basic}) yields
\begin{equation}
\begin{array}{ll}
\di& \di|\tilde \t^{\frac{1}{2}}(x,t)-\tilde\t^{\frac{1}{2}}(b_k(t),t)| \leq C|\tilde{\t}^{\frac{\beta+1}{2}}(x,t)-\tilde{\t}^{\frac{\beta+1}{2}}(b_k(t),t)|
\leq\int_{k}^{k+1}\tilde\t^{\frac{\beta-1}{2}}|\tilde\t_x|dx\\[3mm]
&\di \quad \leq\int_{k}^{k+1}\left(\frac{\t}{\Theta}\right)^{\frac{\beta-1}{2}}\left(\left|\frac{\zeta_x}{\T}\right|+\left|\frac{\zeta\T_x}{\T^2}\right|\right)dx\\[3mm]
&\di\quad\leq C\left(\int_k^{k+1}\frac{\t^\beta\zeta_x^2}{\t^2v}dx\right)^{\frac{1}{2}}
\left(\int_k^{k+1}\frac{\t v}{\T^{\beta+1}}dx\right)^{\frac{1}{2}}+C(M)\left(\int_k^{k+1}\zeta^2\Theta_x^2dx\right)^{\frac{1}{2}}\\[3mm]
&\di\quad\leq C\left(\int_k^{k+1}\frac{\t^\beta\zeta_x^2}{\t^2v}dx\right)^{\frac{1}{2}}
\left(\int_k^{k+1}\t dx\right)^{\frac{1}{2}}\max_{x\in[k,k+1]}v^{\frac{1}{2}}(\cdot,t)
+C(M)\left(\int_k^{k+1}\zeta^2\Theta_x^2dx\right)^{\frac{1}{2}}\\[3mm]
&\di\quad\leq CV^{\frac{1}{2}}(t)\max_{x\in[k,k+1]}v^{\frac{1}{2}}(\cdot,t)
+C(M)\left(\int\zeta^2\Theta_x^2dx\right)^{\frac{1}{2}} 
\end{array}
\end{equation}
for $x\in[k,k+1]$, $k=0$, $\pm1$, $\pm2, \cdots$, 
which along with \eqref{root2} leads to
\begin{equation}\label{v-theta}
\begin{array}{ll}
\di \frac{\a_1}{8}-CV(t)\max_{x\in\mathbb{R}}v(\cdot,t)
-C(M)\int_{\mathbb{R}}\zeta^2\Theta^2_xdx\leq \theta(x,t)\\[3mm]
\di\quad\leq C_0 +CV(t)\max_{x\in\mathbb{R}}v(\cdot,t)
+C(M)\int_{\mathbb{R}}\zeta^2\Theta^2_xdx, 
\end{array}
\end{equation}
for any $x\in\mathbb{R}$. 
Putting this into \eqref{v-repre2} and applying Gronwall's inequality and \eqref{basic} shows that for any $t\in [0,+\infty)$,
\begin{equation}\label{v-upper}
\di v(x,t)\leq C_0.
\end{equation}
Then, integrating \eqref{v-repre} in $x$ over $[k,k+1]$, we obtain after using (\ref{root}) that 
\begin{equation}
\begin{array}{ll}
\di v_-\a_1&\di\leq C_0e^{-t/C_0}+C_0\int_0^t\frac{Y(t)}{Y(s)}\int_k^{k+1}\theta(x,s)dxds\\
&\di \leq C_0e^{-t/C_0}+C_0\int_0^t\frac{Y(t)}{Y(s)}ds,
\end{array}
\end{equation}
which directly yields that
\begin{equation}\label{Yt-Ys}
\di \int_0^t\frac{Y(t)}{Y(s)} ds\geq C_0-C_0e^{-t/C_0}.
\end{equation}
From \eqref{basic},\eqref{v-repre}, \eqref{v-theta}, \eqref{v-upper} and \eqref{Yt-Ys}, and choosing $\d$
suitable small, we have
\begin{equation}\label{wt-1}
\begin{array}{ll}
\di v(x,t)&\di\geq C_0\int_0^t\frac{Y(t)}{Y(s)}\theta(x,s)ds\\
&\di\geq C_0\int_0^t\frac{Y(t)}{Y(s)}ds-C_0\int_0^t\frac{Y(t)}{Y(s)}V(s) dxds-C(C_0,M)\int_0^t\int\zeta^2\Theta_x^2 dxds\\
&\di\geq C_0-C_1(C_0)e^{-t/C_0}-C_0e^{-t/2C_0}\int_0^{t/2}V(s) dxds
 -C_0\int_{t/2}^t V(s)dxds\\
&\di\geq C_0/2
\end{array}
\end{equation}
for any $x\in\mathbb{R}$, provided $t\geq T_0$. Here $C_1(C_0)$ is some positive constant depending on $C_0$, while $T_0$, $C_0$ are positive constants independent of $t$.

Finally, it follows from \cite{kazhikhov} that for any $x\in[k,k+1]$, $k=0,\pm1,\pm2,\cdots$, 
\begin{equation}\label{v-ka}
\di v(x,t)=\frac{1}{\widetilde Y(t)\widetilde B(x,t)}
\left(v_0(x)+\frac{R}{\tilde{\mu}}\int_0^t\widetilde Y(s)\widetilde B(x,s)\t(x,s) ds\right),
\end{equation}
where
$$
\di\widetilde Y(t)=\frac{v_0(a_k(t))}{v(a_k(t),t)}\exp\left(\frac{R}{\tilde{\mu}}\int_0^t\frac{\t}{v}(a_k(t),s)ds\right)
$$
and
$$
\di\widetilde B(x,t)=\exp\left(\frac{1}{\tilde{\mu}}\int_{a_k(t)}^x\Big(u_0(y)-u(y,t)\Big)dy\right)
$$
with $a_{k}(t)$ is the same as in \eqref{root2}.
It follows from \eqref{v-ka} that
\begin{equation}\label{v-ka2}
\di \widetilde Y(t)v(x,t)=\frac{1}{\widetilde B(x,t)}
\left(v_0(x)+\frac{R}{\tilde{\mu}}\int_0^t\widetilde Y(s)\widetilde B(x,s)\t(x,s) ds\right).
\end{equation}
Integrating \eqref{v-ka2} over $[k,k+1]$ with respect to $x$, we obtain
\begin{equation}
\begin{array}{ll}
\di \a_1v_- \widetilde Y(t)
\leq C_0+C_0\int_0^t\widetilde  Y(s)ds,
\end{array}
\end{equation}
which together with the Gronwall's inequality shows that
$$
\widetilde  Y(t)\leq C(C_0,T).
$$
From \eqref{v-ka}, one has
\begin{equation*}
\di v(x,t)\geq \frac{v_0(x)}{\widetilde Y(t)\widetilde B(x,t)}\geq C(C_0,T),
\qquad \forall x\in \mathbb{R},
\end{equation*}
which together with \eqref{wt-1} and \eqref{v-upper}, completes the proof of Lemma \ref{v-bound}.
\hfill $\Box$

\

To obtain the uniform bound of $\t$ from below and above with respect to space and time, motivated by \cite{L-L} and \cite{Li-Xu},  we state the following two lemmas which show the $L^2$-norm (in both space and time) bound of $\psi_x$ and $\t^{-1/2}/\zeta_x$. 

\begin{lemma}\label{lemma-ll-1}
For $\beta\in (0,1)$, there exists some positive constants $C_0$ such that for any given $T>0$,
\begin{equation}\label{new-nergy}
\di\sup_{0\leq t\leq T}\int(\zeta^2+\psi^4)dx+\int_0^T\int((1+\t+\psi^2)\psi_x^2+(\t^{-1}+\t^\beta)\zeta_x^2) dxdt
\leq C_0.
\end{equation}
\end{lemma}
\textbf{Proof}: First, for $t\geq0$, and $a>1$, denoting
\begin{equation}
	\label{omega_a}
	\Omega_a(t)\triangleq \left\{x\in\mathbb{R}\Bigg|\frac{\t}{\Theta}(x,t)>a\right\}
	=\{x\in\mathbb{R}|\zeta(x,t)>(a-1)\Theta(x,t)\}.
\end{equation}

We derive from \eqref{basic} that $\Omega_a$ is bounded since
\begin{equation}\label{omega-bound}
	\di a|\Omega_a|<\sup_{0\leq t\leq T}\int_{\Omega_a}\frac{\t}{\Theta}dx
	\leq C(a)\sup_{0\leq t\leq T}\int\Phi\left(\frac{\t}{\Theta}\right)dx\leq C(a,C_0),
\end{equation}
and that
\begin{equation}
	\di\sup_{0\leq t\leq T}\left(\left|\left(\frac{\t}{\T}>a\right) (t) \right|+\left|\left(\frac{\t}{\T}<a^{-1}\right) (t) \right|  \right) \leq C(a,C_0).
	\label{cedu}
\end{equation}

Next, we divide this proof into three steeps.

Step 1.\quad Multiplying $\eqref{perturb}_3$ by $(\zeta-\Theta)_+=\max\{\zeta-\Theta,0\}$, then integrating the resulted equation
over $\mathbb{R}\times[0,t]$, one has
\begin{equation}\label{5.1}
\begin{array}{ll}
\di\frac{c_{\nu}}{2}\int(\zeta-\Theta)_+^2dx+\tilde{\k}\int_0^t\int_{\Omega_2}\frac{\t^\beta\zeta_x^2}{v}dxds
=\frac{c_{\nu}}{2}\int(\zeta_0(x)-\Theta(x,0))_+^2dx\\[2mm]
\di\quad-\int_0^t\int\frac{R\zeta+R\Theta}{v}\psi_x(\zeta-\Theta)_+ dxds
-\int_0^t\int\frac{R\zeta-p_+\phi}{v}U_x(\zeta-\Theta)_+ dxds\\[2mm]
\di\quad+\tilde{\k}\int_0^t\int_{\Omega_2}\frac{\T^\beta\zeta_x\Theta_x}{V}dxds
-\tilde{\k}\int_0^t\int_{\Omega_2}\frac{(\T^\beta v - \t^\beta V)\Theta_x^2}{vV}dxds
+\tilde{\mu}\int_0^t\int\frac{\psi_x^2}{v}(\zeta-\Theta)_+ dxds\\[2mm]
\di\quad+2\tilde{\mu}\int_0^t\int\frac{\psi_xU_x}{v}(\zeta-\Theta)_+ dxds-\tilde{\mu}\int_0^t\int\frac{\phi U_x^2}{vV}(\zeta-\Theta)_+ dxds\\[2mm]
\di\quad-\int_0^t\int\widetilde R_2(\zeta-\Theta)_+ dxds-c_{\nu}\int_0^t\int\partial_t\Theta(\zeta-\Theta)_+ dxds.
\end{array}
\end{equation}

To estimate the right hand side of (\ref{5.1}), We multiply $\eqref{perturb}_2$ by $2\psi(\zeta-\Theta)_+$, and integrate the resulting equation over
$\mathbb{R}\times[0,t]$ to get
\begin{equation}\label{5.2}
\begin{array}{ll}
\di \int\psi^2(\zeta-\Theta)_+dx+2\tilde{\mu}\int_0^t\int\frac{\psi_x^2}{v}(\zeta-\Theta)_+ dxds
=\int\psi_0^2(x)(\zeta_0(x)-\Theta(x,0))_+dx\\[2mm]
\di\quad+2\int_0^t\int\frac{R\zeta-p_+\phi}{v}\psi_x(\zeta-\Theta)_+dxds
+2\int_0^t\int_{\Omega_2}\frac{R\zeta-p_+\phi}{v}\psi\zeta_xdxds\\[2mm]
\di\quad-2\int_0^t\int_{\Omega_2}\frac{R\zeta-p_+\phi}{v}\psi\Theta_xdxds
+2\tilde{\mu}\int_0^t\int\frac{\phi U_x}{vV}\psi_x(\zeta-\Theta)_+ dxds\\[2mm]
\di\quad-2\tilde{\mu}\int_0^t\int_{\Omega_2}\frac{\psi\psi_x\zeta_x}{v}dxds
+2\tilde{\mu}\int_0^t\int_{\Omega_2}\frac{\phi\psi U_x}{vV}\zeta_xdxds\\[2mm]
\di\quad+2\tilde{\mu}\int_0^t\int_{\Omega_2}\frac{\psi\psi_x\Theta_x}{v}dxds
-2\tilde{\mu}\int_0^t\int_{\Omega_2}\frac{\phi\psi}{vV}U_x\Theta_xdxds\\[2mm]
\di\quad-2\int_0^t\int\psi\widetilde R_1(\zeta-\Theta)_+ dxds
+\int_0^t\int_{\Omega_2}\psi^2\partial_t\zeta dxds-\int_0^t\int_{\Omega_2}\psi^2\partial_t\Theta dxds.
\end{array}
\end{equation}
Adding \eqref{5.2} into \eqref{5.1}, we obtain after using  $\eqref{perturb}_3$ that
\begin{equation}
\begin{array}{ll}\label{long}
\di\int\left(\frac{c_{\nu}}{2}(\zeta-\Theta)_+^2+\psi^2(\zeta-\Theta)_+\right)dx
+\tilde{\mu}\int_0^t\int\frac{\psi_x^2}{v}(\zeta-\Theta)_+ dxds+  \tilde{\kappa}\int_0^t\int_{\Omega_2}\frac{\t^\beta\zeta_x^2}{v}dxds\\[3mm]
\di=\int\left(\frac{c_{\nu}}{2}(\zeta_0(x)-\Theta(x,0))_+^2+\psi_0^2(x)(\zeta_0(x)-\Theta(x,0))_+\right)dx
+\int_0^t\int\frac{R\zeta-2p_+\phi-R\Theta}{v}\psi_x(\zeta-\Theta)_+ dxds\\[3mm]
\di-\int_0^t\int \frac{R\zeta-p_+\phi}{v}U_x(\zeta-\Theta)_+dxds
+\tilde{\k}\int_0^t\int_{\Omega_2}\frac{\T^\beta\zeta_x\Theta_x}{V}dxds
-\tilde{\k}\int_0^t\int_{\Omega_2}\frac{(\T^\beta v - \t^\beta V)\Theta_x^2}{vV}dxds
\\[2mm]
\di+2\tilde{\mu}\int_0^t\int\frac{\psi_xU_x}{V}(\zeta-\Theta)_+dxds
-\tilde{\mu}\int_0^t\int\frac{\phi U_ x^2}{vV}(\zeta-\Theta)_+dxds
+2\int_0^t\int_{\Omega_2}\frac{R\zeta-p_+\phi}{v}\psi\zeta_xdxds\\[3mm]
\di-2\int_0^t\int_{\Omega_2}\frac{R\zeta-p_+\phi}{v}\psi\Theta_x dxds
-2\tilde{\mu}\int_0^t\int_{\Omega_2}\frac{\psi\psi_x\zeta_x}{v}dxds
+2\tilde{\mu}\int_0^t\int_{\Omega_2}\frac{\phi\psi U_x}{vV}\zeta_xdxds\\[3mm]
\di+2\tilde{\mu}\int_0^t\int_{\Omega_2}\frac{\psi\psi_x\Theta_x}{v}dxds
-2\tilde{\mu}\int_0^t\int_{\Omega_2}\frac{\phi\psi}{vV}U_x\Theta_xdxds
-2\int_0^t\int\psi\widetilde R_1(\zeta-\Theta)_+dxds\\[3mm]
\di-\int_0^t\int\widetilde R_2(\zeta-\Theta)_+dxds
-c_{\nu}\int_0^t\int\partial_t\Theta(\zeta-\Theta)_+dxds
-\int_0^t\int_{\Omega_2}\psi^2\partial_t\Theta dxds\\[3mm]
\di+\frac{\tilde{\mu}}{c_{\nu}}\int_0^t\int_{\Omega_2}\psi^2\left(\frac{\psi_x^2+2\psi_xU_x}{v}-\frac{\phi U_x^2}{vV}\right)dxds
-\frac{1}{c_{\nu}}\int_0^t\int_{\Omega_2}\psi^2\widetilde R_2dxds\\[3mm]
\di-\frac{1}{c_{\nu}}\int_0^t\int_{\Omega_2}\psi^2\left (\frac{R\zeta+R\Theta}{v}\psi_x+\frac{R\zeta-p_+\phi}{v}U_x\right)dxds
+\frac{\tilde{\k}}{c_{\nu}}\int_0^t\int_{\Omega_2}\psi^2\left(\frac{\t^\beta\t_x}{v}-\frac{\T^\beta\Theta_x}{V}\right)_xdxds\\[3mm]
\di\triangleq \int\left(\frac{c_{\nu}}{2}(\zeta_0(x)-\Theta(x,0))_+^2+\psi_0^2(x)(\zeta_0(x)-\Theta(x,0))_+\right)dx
+\sum_{i=1}^{20}I_i.
\end{array}
\end{equation}
We will estimate \eqref{long} term by term. Recalling \eqref{basic}, \eqref{v} and \eqref{omega-bound}, it holds that
\begin{equation}
\begin{array}{ll}\label{i1}
\di|I_1|
&\di=\left|\int_0^t\int\frac{R\zeta-2p_+\phi-R\Theta}{v}\psi_x(\zeta-\Theta)_+ dxds\right|\\[2mm]
&\di\leq \frac{\tilde{\mu}}{4}\int_0^t\int\frac{\psi_x^2}{v}(\zeta-\Theta)_+ dxds
+C_0\int_0^t\int(\zeta^2+\phi^2\zeta)(\zeta-\Theta)_+dxds\\
&\di\leq \frac{\tilde{\mu}}{4}\int_0^t\int\frac{\psi_x^2}{v}(\zeta-\Theta)_+ dxds
+C_0\int_0^t\int(\zeta+\phi^2)\left(\zeta-\frac{1}{2}\Theta\right)_+^2dxds\\[2mm]
&\di\leq \frac{\tilde{\mu}}{4}\int_0^t\int\frac{\psi_x^2}{v}(\zeta-\Theta)_+ dxds
+C_0\int_0^t\max_{x\in\mathbb{R}}\left(\zeta-\frac{1}{2}\Theta\right)_+^2\int_{\{\zeta>\frac{\Theta}{2}\}}(\zeta+\phi^2)dxds\\
&\di
\leq \frac{\tilde{\mu}}{4}\int_0^t\int\frac{\psi_x^2}{v}(\zeta-\Theta)_+ dxds
+C_0\int_0^t\max_{x\in\mathbb{R}}\left(\zeta-\frac{1}{2}\Theta\right)_+^2ds.
\end{array}
\end{equation}

It follows from \eqref{important}, \eqref{basic}, \eqref{phi-x3},\eqref{omega_a} and Cauchy's inequality that
\begin{equation}
\begin{array}{ll}\label{i3}
\di|I_2|+|I_3|+|I_4|\\
\di\leq C(M)\int_0^t\int(\phi^2+\zeta^2)|U_x|dxds+
\tilde{\k}\int_0^t\int_{\Omega_2}\left|\frac{\zeta}{\Theta}\right|\left(\left|\frac{\T^\beta\zeta_x\Theta_x}{V}\right|
+\left|\frac{(\T^\beta v-\t^\beta V)\Theta_x^2}{vV}\right|\right) dxds\\[3mm]
\di\leq C(M)\int_0^t\int(\phi^2+\zeta^2)|U_x|dxds + \frac{\tilde{\kappa}}{8} \int_{0}^{t}\int_{\Omega_2}\frac{\t^\beta\zeta_x^2}{v}dxds+C\int_{0}^{t}\int_{\Omega_2}\frac{\T^{2\beta}\T_x^2v}{V^2\t^\beta}dxds\\[3mm]
\di\quad+C(M)\int_{0}^{t}\int_{\Omega_2}\zeta^2\T_x^2dxds\\[ 3mm]
\di\leq \frac{\tilde{\k}}{8}\int_0^t\int_{\Omega_2}\frac{\t^\beta\zeta_x^2}{v}dxds
+C(M)\int_0^t\int(\phi^2+\zeta^2)(|U_x|+\Theta_x^2)dxds\\
\di\leq \frac{\tilde{\k}}{8}\int_0^t\int_{\Omega_2}\frac{\t^\beta\zeta_x^2}{v}dxds+C_0.
\end{array}
\end{equation}
Similarly, one has
\begin{equation}
\begin{array}{ll}\label{i4}
\di|I_5|+|I_6|+|I_{10}|+|I_{11}|+|I_{12}|\\
\di\leq C_0\int_0^t\int\frac{\psi_x^2}{\t}dxds
+\frac{\tilde{\k}}{8}\int_0^t\int_{\Omega_2}\frac{\t^\beta\zeta_x^2}{v}dxds\\
\di\quad+C(C_0,M)\int_0^t\int(\phi^2+\psi^2+\zeta^2)(U_x^2+\Theta_x^2)dxds\\
\di\leq\frac{\tilde{\k}}{8}\int_0^t\int_{\Omega_2}\frac{\t^\beta\zeta_x^2}{v}dxds+C_0.
\end{array}
\end{equation}
It holds that By Cauchy's inequality, \eqref{basic}, \eqref{v} and \eqref{omega-bound}, 
\begin{equation}
\begin{array}{ll}\label{i7}
\di|I_7|\leq\frac{\tilde{\k}}{8}\int_0^t\int_{\Omega_2}\frac{\t^\beta\zeta_x^2}{v}dxds
+C_0\int_0^t\int_{\Omega_2}(\zeta^2+\phi^2)\psi^2dxds\\
\di\quad\leq\frac{\tilde{\k}}{8}\int_0^t\int_{\Omega_2}\frac{\t^\beta\zeta_x^2}{v}dxds
+C_0\int_0^t\int_{\Omega_2}(\zeta^2\psi^2+\phi^2\psi^4+\phi^2)dxds\\
\di\quad\leq\frac{\tilde{\k}}{8}\int_0^t\int_{\Omega_2}\frac{\t^\beta\zeta_x^2}{v}dxds
+C_0\int_0^t\int_{\Omega_2}(\zeta^2(\phi^2+\psi^2)+\phi^2\psi^4)dxds\\
\di\quad\leq\frac{\tilde{\k}}{8}\int_0^t\int_{\Omega_2}\frac{\t^\beta\zeta_x^2}{v}dxds
+C_0\int_0^t\int_{\Omega_2}\left(\left(\zeta-\frac{1}{2}\Theta\right)^2(\phi^2+\psi^2)+\phi^2\psi^4\right)dxds\\
\di\quad\leq\frac{\tilde{\k}}{8}\int_0^t\int_{\Omega_2}\frac{\t^\beta\zeta_x^2}{v}dxds
+C_0\int_0^t\left(\max_{x\in\mathbb{R}}\left(\zeta-\frac{1}{2}\Theta\right)_+^2+\max_{x\in\mathbb{R}}\psi^4\right)\int(\phi^2+\psi^2)dxds\\
\di\quad\leq\frac{\tilde{\k}}{8}\int_0^t\int_{\Omega_2}\frac{\t^\beta\zeta_x^2}{v}dxds
+C_0\int_0^t\left(\max_{x\in\mathbb{R}}\left(\zeta-\frac{1}{2}\Theta\right)_+^2+\max_{x\in\mathbb{R}}\psi^4\right)ds.
\end{array}
\end{equation}
Similarly, 
\begin{equation}
\begin{array}{ll}\label{i8}
\di|I_8|\leq C_0\int_0^t\int_{\Omega_2}(\psi^2+(\phi^2+\zeta^2)\Theta_x^2)dxds\\
\di\quad\leq C_0\int_0^t\int_{\Omega_2}(\zeta^2\psi^2+(\phi^2+\zeta^2)\Theta_x^2)dxds\\
\di\quad\leq C_0\int_0^t\max_{x\in\mathbb{R}}\left(\zeta-\frac{1}{2}\Theta\right)_+^2\int\psi^2dxds
+C_0\int_0^t\int(\phi^2+\zeta^2)\Theta_x^2dxds\\
\di\quad\leq C_0\int_0^t\max_{x\in\mathbb{R}}\left(\zeta-\frac{1}{2}\Theta\right)_+^2ds+C_0.
\end{array}
\end{equation}
By Cauchy's inequality, it holds that
\begin{equation}
\begin{array}{ll}\label{i9}
\di|I_9|\leq \frac{\tilde{\k}}{8}\int_0^t\int_{\Omega_2}\frac{\t^\beta\zeta_x^2}{v}dxds
+C_0\int_0^t\int_{\Omega_2}\psi^2\psi_x^2dxds.
\end{array}
\end{equation}
Recalling Lemma \ref{decay}, \eqref{basic} and \eqref{phi-x3}, one has
\begin{equation}
\begin{array}{ll}\label{i10}
\di|I_{13}|+|I_{14}|+|I_{15}|+|I_{16}|+|I_{18}|\\
\di\leq C\d \int_0^t\int (1+s)^{-1}(\psi^2+\zeta^2)e^{-\frac{c_1x^2}{1+s}}dxds
\leq C_0.
\end{array}
\end{equation}
By Cauchy's inequality, and using \eqref{basic}, \eqref{phi-x3} and Lemma \ref{decay}, we obtain
\begin{equation}
\begin{array}{ll}\label{i17}
\di|I_{17}|\leq C_0\int_0^t\int(\psi^2\psi_x^2+\psi^2|\psi_x||U_x|+\psi^2|\phi|U_x^2)dxds\\
\di\quad\leq C_0\int_0^t\int(\psi^2\psi_x^2+\psi^2U_x^2+\psi^2|\phi|U_x^2)dxds\\
\di\quad\leq C_0\int_0^t\int\psi^2\psi_x^2dxds+C(C_0,M)\int_0^t\int\psi^2U_x^2dxds\\
\di\quad\leq C_0\int_0^t\int\psi^2\psi_x^2dxds+C_0.
\end{array}
\end{equation}
Similarly,
\begin{equation}
\begin{array}{ll}\label{i19}
\di|I_{19}|\leq\int_0^t\int_{\Omega_2}(\psi^2\psi_x^2+\psi^2\zeta^2+\psi^4|U_x|+(\zeta^2+\phi^2)|U_x|)dxds\\
\di\leq C_0\int_0^t\int\psi^2\psi_x^2dxds
+C_0\int_0^t\max_{x\in\mathbb{R}}\left(\zeta-\frac{1}{2}\Theta\right)_+^2\int\psi^2dxds\\
\di\quad+\int_0^t\max_{x\in\mathbb{R}}\psi^4\int|U_x|dxds+C_0\\
\di\leq C_0\int_0^t\int\psi^2\psi_x^2dxds
+C_0\int_0^t\left(\max_{x\in\mathbb{R}}\left(\zeta-\frac{1}{2}\Theta\right)_+^2+\max_{x\in\mathbb{R}}\psi^4\right)ds
+C_0.
\end{array}
\end{equation}
Finally, for $\eta>0$ and 
\begin{equation}\label{cut-off}
\di \varphi_{\eta}(z)=\left\{
\begin{array}{ll}
1, & z>\eta,\\
z/\eta, & 0< z\leq \eta,\\
0  & z\leq 0.
\end{array}
\right.
\end{equation}
Integrating by parts gives
\begin{equation}
\begin{array}{ll}\label{i20}
\di I_{20}=\frac{\tilde{\kappa}}{c_{\nu}}\int_0^t\int_{\Omega_2}\psi^2\left(\frac{\t^\beta\t_x}{v}-\frac{\T^\beta\Theta_x}{V}\right)_xdxds\\
\di=\frac{\tilde{\k}}{c_{\nu}}\lim_{\eta\rightarrow0+}\int_0^t\int\varphi_{\eta}(\zeta-\Theta)\psi^2\left(\frac{\t^\beta\zeta_x-\t^\beta\Theta_x}{v}\right)_xdxds\\
\di\quad+\frac{\tilde{\k}}{c_{\nu}}\int_0^t\int_{\Omega_2}\psi^2\left(\frac{2\t^\beta\Theta_x}{v}-\frac{\T^\beta\Theta_x}{V}\right)_xdxds
\triangleq I_{20}^1+I_{20}^2.
\end{array}
\end{equation}
Lebegue's dominated convergence theorem shows that for $\beta<1$ and any $\varepsilon>0$
\begin{equation}
\begin{array}{ll}\label{i201}
\di I_{20}^1&\di=-\frac{2\tilde{\kappa}}{c_{\nu}}\lim_{\eta\rightarrow 0+}\int_0^t\int\varphi_{\eta}(\zeta-\Theta)\psi\psi_x\frac{\t^\beta(\zeta_x-\Theta_x)}{v}dxds\\[3mm]
&\di\quad-\frac{\tilde{\k}}{c_{\nu}}\lim_{\eta\rightarrow 0+}\int_0^t\int\varphi_{\eta}'(\zeta-\Theta)\frac{\psi^2\t^\beta(\zeta_x-\Theta_x)^2}{v}dxds\\[3mm]
&\di\leq C\int_0^t\int_{\Omega_2}\frac{|\psi\psi_x\t^\beta\zeta_x|+|\psi\psi_x\t^\beta\Theta_x|}{v} dxds\\[3mm]
&\di\leq C(\varepsilon)\int_{0}^{t}\int|\psi|^{2/(1-\beta)}\psi_x^2dxds + \varepsilon\int_{0}^{t}\int\left(\psi_x^2\t+\t^\beta\zeta^2_x\right)dxds+ \varepsilon\int_{0}^{t}\int_{\Omega_2}\t^\beta\T^2_xdxds,
\end{array}
\end{equation}
where in the second inequality we have used both $\varphi_{\eta}(z)\in [0,1]$ and $\varphi_{\eta}'(z)\geq 0$.


 Noticing that
 \begin{equation}
 	\begin{array}{ll}\label{5.4}
 		\di\int_0^t\int(\t\psi_x^2+\t^\beta\zeta_x^2)dxds=\int_0^t\int_{\{\zeta>2\Theta\}}(\t\psi_x^2+\t^\beta\zeta_x^2)dxds
 		+\int_0^t\int_{\{\zeta\leq 2\Theta\}}(\t\psi_x^2+\t^\beta\zeta_x^2)dxds\\[3mm]
 		\di\quad\leq\int_0^t\int_{\{\zeta>2\Theta\}}\left(\frac{3}{2}\psi_x^2\zeta+C_0\frac{\t^\beta\zeta_x^2}{v}\right) dxds
 		+\int_0^t\int_{\{\zeta\leq 2\Theta\}}\left(\frac{\psi_x^2}{\t}+\frac{\t^\beta\zeta_x^2}{\t^2}\right)\t^2dxds\\[3mm]
 		\di\quad\leq\int_0^t\int_{\{\zeta>2\Theta\}}\left(3\psi_x^2(\zeta-\Theta)+C_0\frac{\t^\beta\zeta_x^2}{v}\right) dxds
 		+C\int_0^t\int_{\{\zeta\leq 2\Theta\}}\left(\frac{\psi_x^2}{\t}+\frac{\t^\beta\zeta_x^2}{\t^2}\right)dxds\\[3mm]
 		\di\quad\leq C_0\int_0^t\int_{\Omega_2}\left(\frac{\psi_x^2}{v}(\zeta-\Theta)_+ +\frac{\t^\beta\zeta_x^2}{v}\right)dxds+C_0,
 	\end{array}
 \end{equation}
 and 
\begin{equation}
	 \int_{0}^{t}\int_{\Omega_2}\t^\beta\T^2_xdxds\leq C(M)\int_{0}^{t}\int_{\Omega_2}\zeta^2\T^2_xdxds\leq C(M)\delta\leq C_0,
\end{equation}
where we choose $\delta$ suitable small. 


Similarly,
\begin{equation}
\begin{array}{ll}\label{i202}
\di I_{20}^2=\frac{\k}{c_{\nu}}\int_0^t\int_{\Omega_2}\psi^2\left(\frac{2({\t^\beta\Theta_{x}})_x}{v}-\frac{(\T^\beta\Theta_{x})_x}{V}
-\frac{2\t^\beta\Theta_xv_x}{v^2}+\frac{\T^\beta\Theta_xV_x}{V^2}\right) dxds\\
\di\quad\leq C_0\int_0^t\int \psi^2\left(\beta\t^{\beta-1}|\t_x||\T_x|+\t^\beta|\T_{xx}|+\beta\T^{\beta-1}|\T_x|+\T^\beta|\T_{xx}|\right)dxds\\[3mm]
\di\quad \quad+C_0\int_0^t\int \psi^2\left(\t^\beta|\T_x||v_x|+\T^\beta|\T_x||V_x|\right)dxds
\\[3mm]
\di\quad\leq C_0\int_0^t\int\frac{\psi^2|\zeta_x||\T_x|}{\t^{1-\beta}}+\frac{\psi^2\T_x^2}{\t^{1-\beta}} dxds+C_0\int_0^t\int\psi^4\left(|\T_{xx}|+|\T_x|+\T_x^2\right)dxds\\[3mm]
\di\quad\quad+C(M)\int_0^t\int_{\Omega_2}\zeta^2\left(|\T_{xx}|+|\T_x| +  V_x^2\right)dxds + C(M)\delta\int_0^t\int\frac{\t\phi_x^2}{v^3}dxds

\\[3mm]

\di\quad\leq \frac{\tilde{\k}}{8}\int_0^t\int\frac{\t^\beta\zeta_x^2}{v}dxds +C_0\int_{0}^{t}\max_{x\in\mathbb{R}}\psi^4\left(\int|\T_{xx}|+|\T_x|+\T_x^2dx\right)ds + C_0\\[3mm]

\di\quad\leq \frac{\tilde{\k}}{8}\int_0^t\int\frac{\t^\beta\zeta_x^2}{v}dxds+C_0\int_{0}^{t}\max_{x\in\mathbb{R}}\psi^4ds +C_0.

%
\end{array}
\end{equation}

Substituting the estimates \eqref{i1}-\eqref{i202} into \eqref{long}, and using \eqref{5.4}, we have
\begin{equation}
\begin{array}{ll}\label{5.3}
\di \int(\zeta-\Theta)_+^2
 dx+\int_0^t\int(\t\psi_x^2+\t^\beta\zeta_x^2)dxds \leq C_0 \\
\di\quad 
+C_0\int_0^t\max_{x\in\mathbb{R}}\left(\zeta-\frac{1}{2}\Theta\right)_+^2ds+  C\int_{0}^{t}\int|\psi|^{2/(1-\beta)}\psi_x^2dxds ,
\end{array}
\end{equation}
where we have used the simple fact by using (\ref{basic}) that for any $\delta>0$
\begin{equation*}
	\begin{array}{ll}
		&\di\int_0^t\int\psi^2\psi_x^2dxds\leq\int_{0}^{t}\int\psi_x^2dxds + \int_{0}^{t}\int|\psi|^{2/(1-\beta)}\psi_x^2dxds\\[3mm]
		&\di\quad\leq\delta\int_{0}^{t}\int\t\psi_x^2dxds + C(\delta)\int_{0}^{t}\int\frac{\psi^2_x}{\t v}dxds+ \int_{0}^{t}\int|\psi|^{2/(1-\beta)}\psi_x^2dxds,
	\end{array}
\end{equation*}
and
\begin{equation*}
	\begin{array}{ll}
		\di\int_{0}^{t}\max_{x\in\mathbb{R}}\psi^4ds = \int_{0}^{t}\left(\max_{x\in\mathbb{R}}\int_{x}^{\infty}\partial_y\psi^2dy\right)^2ds\\[3mm]
		\di\quad\leq C\int_{0}^{t}\left(\int|\psi\psi_x|dx\right)^2ds\leq \int_{0}^{t}\|\psi\|^2_{L^2}\|\psi_x\|^2_{L^2}ds\leq C\int_{0}^{t}\|\psi_x\|^2_{L^2}ds.
	\end{array}
\end{equation*}

\

\underline{ Step 2. }\quad To estimate the last term on the right hand side of \eqref{5.3},
we multiply $\eqref{perturb}_2$ by $|\psi|^b\psi(b=2/(1-\beta))$ and integrate the resulting equality  over $\mathbb{R}\times[0,t]$ to get
\begin{equation}
\begin{array}{ll}\label{5.5}
\di\frac{1}{b+2}\int\psi^{b+2}dx+(b+1)\tilde{\mu}\int_0^t\int\frac{|\psi|^b\psi_x^2}{v} dxds
=\frac{1}{b+2}\int\psi_0^{b+2}dx+(b+1)R\int_0^t\int\frac{\zeta|\psi|^b\psi_x}{v}dxds\\[3mm]
\di\quad\quad-(b+1))p_+\int_0^t\int\frac{\phi|\psi|^b\psi_x}{v}dxds+(b+1)\tilde{\mu}\int_0^t\int\frac{\phi U_x}{vV}|\psi|^b\psi_x dxds
-\int_0^t\int\widetilde R_1|\psi|^b\psi dxds\\
\di\quad\leq C+C \int_{0}^{t}\int\left(\frac{|\zeta|}{v}1_{(\zeta<3\T)}+\frac{|\phi|}{v}\right)|\psi|^b|\psi_x|dxds + C \int_{0}^{t}\int_{\Omega_3}\frac{\zeta}{v}|\psi|^b|\psi_x|dxds\\[3mm]
\di\quad\quad+C\int_0^t\int\frac{\phi U_x}{vV}|\psi|^b\psi_x dxds
-\int_0^t\int\widetilde R_1|\psi|^b\psi dxds\\[3mm]
\di\quad\triangleq C_0+\sum_{i=1}^4J_i.
\end{array}
\end{equation}

It follows from (\ref{basic}) and (\ref{v}) that for any $\delta_1 \in [2\T,4\Theta]$,
\begin{equation}
	\begin{array}{ll}
		\di\sup_{0\leq t\leq T}\int\phi^2dx + \sup_{0\leq t\leq T}\int_{\{\zeta<\delta_1\}}\zeta^2dx\\[3mm]
		\di\quad \leq C\sup_{0\leq t\leq T}\int(\widetilde{v}-1)^2dx+C\sup_{0\leq t\leq T}\int_{\{\zeta<\delta_1\}}(\widetilde{\t}-1)^2dx\\[3mm]
		\di\quad \leq C\sup_{0\leq t\leq T}\int(\widetilde{v}-\ln\widetilde{v}-1)dx +  C\sup_{0\leq t\leq T}\int(\widetilde{\t}-\ln\widetilde{\t}-1)dx\leq C_0,
		\label{LL2}
	\end{array}
\end{equation}
which together with Holder's inequality yields that
\begin{equation}
	\begin{array}{ll}
		\label{J11}
	\di|J_1|&\di \leq C\int_{0}^{t}\||\psi|^b \psi_x\|_{L^2}\left(\int\phi^2dx + \int_{\{\zeta<3\T\}}\zeta^2dx\right)^{1/2}ds\\[3mm]
	&\di\leq C\int_{0}^{t}\||\psi|^b \psi_x\|_{L^2}ds.
	\end{array}
\end{equation}

Nest, on the one hand, if $\beta\leq1/2$,
\begin{equation}
	\begin{array}{ll}
		\label{J12}
	\di	\||\psi|^b \psi_x\|_{L^2}&\di\leq C\max_{x\in\mathbb{R}}\psi^2(x,t)\||\psi|^{b\beta} \psi_x\|_{L^2}\\[3mm]&\di\leq C(\varepsilon)\int\psi_x^2dx + \varepsilon\int|\psi|^b \psi_x^2dx,
	\end{array}
\end{equation}
where in the second inequality we hace used (\ref{basic}) and the following simple fact that for any $w\in H^1$,
\begin{equation}
	\begin{array}{ll}
		\label{J13}
		\di \max_{x\in\mathbb{R}} w^2(x) = \max_{x\in\mathbb{R}}\left(-2\int_{x}^{\infty}w(y)w_x(y)dy\right)\leq2\|w\|_{L^2}\|w_x\|_{L^2}.
	\end{array}
\end{equation}
On the other hand, if $\beta\in(1/2,1)$, we have 
\begin{equation}
	\begin{array}{ll}
		\label{J14}
		\di\||\psi|^b\psi_x\|_{L^2}&\di\leq C(\varepsilon)\max_{x\in\mathbb{R}}|\psi|^b(x,t)+\varepsilon\int|\psi|^b \psi_x^2dx\\[3mm]
		&\di\leq C(\varepsilon)\int \psi_x^2dx+2\varepsilon\int|\psi|^b\psi_x^2dx,
	\end{array}
\end{equation}
where we have used 
\begin{equation*}
	\begin{array}{ll}
		\di\max_{x\in\mathbb{R}}|\psi|^b&\di=\max_{x\in\mathbb{R}}\int_{x}^{\infty}\left(-\partial_x|\psi|^b\right)dx\\[3mm]
		&\di\leq C\int|\psi_x||\psi|^{b-1}dx\\[3mm]
		&\di\leq C\left(\int|\psi|^{b-4}\psi_x^2dx\right)^{1/2}\left(\int|\psi|^{b+2}dx\right)^{1/2}\\[3mm]
		&\di\leq C\left(\int|\psi|^{b-4}\psi_x^2dx\right)^{1/2}\max_{x\in\mathbb{R}}|\psi|^{b/2},
	\end{array}
\end{equation*}
due to $b = 2/(1-\beta>4)$. Thus, combining (\ref{J11}), (\ref{J12}) and (\ref{J14}) implies for $\beta\in (0,1)$,
\begin{equation}
	\begin{array}{ll}
		\label{J1-C}
		\di|J_1|\leq C(\varepsilon)\int_{0}^{t}\int\psi_x^2dxds + C\varepsilon\int_{0}^{t}\int|\psi|^b\psi_x^2dxds.
			\end{array}
\end{equation}

Next, if follows from Cauchy's inequality that
\begin{equation}
	\begin{array}{ll}
		\di|J_2|&\di\leq C(\varepsilon)\int_{0}^{t}\int_{\Omega_3}\zeta^2|\psi|^b dxds+\varepsilon\int_{0}^{t}\int\frac{|\psi|^b\psi^2_x}{v}dxds\\[3mm]
		&\di\leq C(\varepsilon)\int_{0}^{t}\max_{x\in\mathbb{R}}(\zeta-\frac{3}{4}\T)_+^{\frac{2b}{b+2}}\int\left(\zeta^2+|\psi|^{b+2}\right)dxds +\varepsilon \int_{0}^{t}\int\frac{|\psi|^b\psi^2_x}{v}dxds\\[3mm]
		&\di\leq C(\varepsilon)\int_{0}^{t}\max_{x\in\mathbb{R}}\left((\zeta-\frac{1}{2}\T)_+^{\beta+2}\zeta^{-1}\right)\int\left(\zeta^2+|\psi|^{b+2}\right)dxds+\varepsilon\int_{0}^{t}\int\frac{|\psi|^b\psi^2_x}{v}dxds,\\[3mm]
		\label{J21}
	\end{array}
\end{equation}
due to $\beta\in (0,1)$ and  $2b/(b+2)<\beta+1$.

Then, directing computation gives that
\begin{equation*}
	\begin{array}{ll}
		\di\max_{x\in\mathbb{R}}\left((\zeta-\frac{1}{2}\T)_+^{\beta+2}\zeta^{-1}\right) = \max_{x\in\mathbb{R}}\int_{x}^{\infty}\left(-\partial_x\left((\zeta-\frac{1}{2}\T)_+^{\beta+2}\zeta^{-1}\right)\right)dx\\[3mm]
		\di\quad\leq C\int\left(|\zeta_x|+|\frac{1}{2}\T_x|\right)(\zeta-\frac{1}{2}\T)_+^{\beta+1}\zeta^{-1} dx+C\int|\zeta_x|(\zeta-\frac{1}{2}\T)_+^{\beta+2}\zeta^{-2}dx\\[3mm]
		\di\quad\leq C\left(\int_{\{\zeta>1/2\T\}}\frac{\t^{2-\beta}}{\zeta^{2-\beta}}\t^{\beta-2}\zeta_x^2dx\right)^{1/2}\left(\int(\zeta-\frac{1}{2}\T)_+^{\beta+2}dx\right)^{1/2}\\[3mm]
		\di\quad\quad+C\left(\int_{\{\zeta>1/2\T\}}\zeta^{\beta-2}\T_x^2dx\right)^{1/2}\left(\int (\zeta-\frac{1}{2}\T)_+^{\beta+2}dx\right)^{1/2}\\[3mm]
		\di\quad\leq C\left(\int\t^{\beta-2}\zeta_x^2dx\right)^{1/2}\max_{x\in\mathbb{R}}\left((\zeta-\frac{1}{2}\T)_+^{\beta+2}\zeta^{-1}\right)^{1/2}\left(\int_{\{\zeta>1/2\T\}}\zeta dx\right)^{1/2}\\[3mm]
		\di\quad\quad + C \left(\int\zeta^{2}\T_x^2dx\right)^{1/2}\max_{x\in\mathbb{R}}\left((\zeta-\frac{1}{2}\T)_+^{\beta+2}\zeta^{-1}\right)^{1/2}\left(\int_{\{\zeta>1/2\T\}}\zeta dx\right)^{1/2}\\[3mm]
		\di\quad\leq C\int\t^{\beta-2}\zeta_x^2dx+ C \int\zeta^{2}\T_x^2dx+\frac{1}{2}\max_{x\in\mathbb{R}}\left((\zeta-\frac{1}{2}\T)_+^{\beta+2}\zeta^{-1}\right),\\[3mm]
	\end{array}
\end{equation*}
where in the last inequality we have used Young's inequality. This implies 
\begin{equation}
	\begin{array}{ll}
		\di\max_{x\in\mathbb{R}}\left((\zeta-\frac{1}{2}\T)_+^{\beta+2}\zeta^{-1}\right)\leq C\int\t^{\beta-2}\zeta_x^2dx+C\int\zeta^2\T_x^2dx.
		\label{J22}
	\end{array}
\end{equation}

From \eqref{basic}, \eqref{important} and \eqref{phi-x3}, we have
\begin{equation}
\begin{array}{ll}\label{j3}
\di|J_3|&\di\leq \tilde{\mu}\int_0^t\int\frac{|\psi|^{b}\psi_x^2}{v}dxds
+C_0\int_0^t\int\phi^2|\psi|^{b}U_x^2dxds\\
&\di\leq\tilde{\mu}\int_0^t\int\frac{|\psi|^{b}\psi_x^2}{v}dxds
+C(C_0,M)\int_0^t\int\phi^2U_x^2dxds\\
&\di\leq\tilde{\mu}\int_0^t\int\frac{|\psi|^{b}\psi_x^2}{v}dxds+C_0.
\end{array}
\end{equation}
Then, combining Lemma \ref{decay}, \eqref{basic}, \eqref{important} and \eqref{phi-x3}, we have
\begin{equation}
\begin{array}{ll}\label{j4}
\di|J_4|&\di\leq O(1)\d\int_0^t\int(1+s)^{-1}|\psi|^{b+1}e^{-\frac{c_1x^2}{1+s}}dxds\\
&\di\leq C(M)\d\int_0^t\int(1+s)^{-1}\psi^2e^{-\frac{c_1x^2}{1+s}}dxds\leq C_0.
\end{array}
\end{equation}
Putting the estimates \eqref{J1-C}-\eqref{j4} into \eqref{5.5} and choose $\varepsilon$ suitable small gives
\begin{equation}
	\begin{array}{ll}
		\di\sup_{0\leq t\leq T}\int|\psi|^{b+2}dx + \int_{0}^{t}\int|\psi|^b\psi_x^2dxds\\[3mm]
		\di \quad\leq C_0+C\int_{0}^{t}\left(\int\t^{\beta-2}\zeta_x^2dx\int\left(\zeta^2+|\psi|^{b+2}\right)dx\right)ds+C\int_{0}^{t}\int\psi_x^2dxds\\[3mm]
		\di\quad\leq C_0\int_{0}^{t}\left(\int\t^{\beta-2}\zeta_x^2dx\int\left(\zeta^2+|\psi|^{b+2}\right)dx\right)ds  + C_0\delta\int_{0}^{t}\int\t\psi_x^2dxds+ C_0(\delta),
	\end{array}
	\label{s2}
\end{equation}
where in the last inequality we have used Young's inequality and (\ref{basic}). 

Adding (\ref{s2}) multiplies by $C_2+1$ to (\ref{5.3}), then choosing $\delta$ suitable small, we have 
\begin{equation}
	\begin{array}{ll}
		\di\sup_{0\leq t\leq T}\int[\zeta^2+|\psi|^{b+2}]dx + \int_{0}^{t}\int[(\t+|\psi|^b)\psi_x^2+\t^\beta\zeta_x^2]dxds\\[3mm]
		\di\leq C_0+C_0\int_{0}^{t}\left(\int\t^{\beta-2}\zeta_x^2dx\int\left(\zeta^2+|\psi|^{b+2}\right)dx\right)ds + C_0\int_{0}^{t}\max_{x\in\mathbb{R}}\left(\zeta-\frac{1}{2}\T\right)_+^2ds.
		\label{5.7}
	\end{array}
\end{equation}
where we have used 
\begin{equation*}
	\begin{array}{ll}
		\di\int\zeta^2dx\leq C\int_{\Omega_3}\zeta^2dx + C\int_{\{\zeta<3\T\}}\zeta^2dx\leq C\int(\zeta-\T)_+^2 + C_0,
	\end{array}
\end{equation*}
due to (\ref{LL2}).

\

\underline{ Step 3. }\quad It remains to estimate the last terms on the right hand side of \eqref{5.7}. Indeed, chooseing $\chi=-1$ in (\ref{psx22})yields that for any $\varepsilon>0$
\begin{equation}
	\begin{array}{ll}
	\di \int_{0}^{t}\max_{x\in\mathbb{R}}(\zeta-\frac{1}{2}\T)_+^2ds &\di=\int_{0}^{t}\max_{x\in\mathbb{R}}(\t-\frac{3}{2}\T)_+^2ds\\[3mm]
	 &\di \leq C(\varepsilon)\int_{0}^{t}\int\frac{\t^\beta\zeta_x^2}{v\t^2}dxds + \varepsilon\int_{0}^{t}\int\t^\beta\zeta_x^2dxds + C_0 \\[3mm]
	&\di \leq C_0 + \varepsilon\int_{0}^{t}\int\t^\beta\zeta_x^2dxds
	\end{array}
\end{equation}
due to $\beta<1$, (\ref{cedu}) and (\ref{basic}). Putting this into (\ref{5.7}), choosing $\varepsilon$ suitable small, and using Gronwall's inequality lead to 
\begin{equation}
	\di\sup_{0\leq t\leq T}\int[\zeta^2+|\psi|^{b+2}]dx + \int_{0}^{t}\int[(\t+|\psi|^b)\psi_x^2+\t^\beta\zeta_x^2]dxds\leq C_0.
\end{equation}
Thus combining this, (\ref{5.3}) and (\ref{basic}) immediately gives(\ref{new-nergy}). The proof of Lemma \ref{lemma-ll-1} is completed.

\
\begin{lemma}\label{beta>1}
For $\beta\in[1,\infty)$, there exsits some positive constant $C$ such that for any $T>0$,
\begin{equation}
	\begin{array}{ll}
		\label{beta-new}
\di \int_{0}^{T}\int(\psi_x^2 + \t^{-1}\zeta_x^2)dxdt\leq C_0.
	\end{array}
\end{equation}
\end{lemma}

\textbf{Proof}:\quad First, using $(\ref{ns})_1$, $(\ref{ns})_2$ and  (\ref{p-e}), we write  $(\ref{ns})_3$ as 
\begin{equation}
	\label{ener}
	c_v\t_t + \frac{R\t}{v}u_x = \tilde{\kappa}\left(\frac{\theta^\beta\t_x}{v}\right)_x +\tilde{\mu} \frac{u_x^2}{v}.
\end{equation}
For $q\geq\beta+4$, multiplying (\ref{ener}) by $(\t^{-q}-4\T^{-q})_+=\max\{(\t^{-q}-4\T^{-q}),0\}$ and integrating the resulting equality over $\mathbb{R} \times (0,t)$ gives
 \begin{equation}
 	\begin{array}{ll}
 		\label{beta-1}
 		\di c_v\int\int_{\t}^{4^{-1/q}\T}(y^{-q}-4\T^{-q})_+dydx + q\tilde{\kappa}\int_{0}^{t}\int_{\{\t<4^{-1/q}\T\}}v^{-1}\t^{\beta-q-1}\t_x^2dxds\\[3mm]
 		\di\quad +\tilde{\mu} \int_{0}^{t}\int v^{-1}u_x^2(\t^{-q}-4\T^{-q})_+dxds\\[3mm]
 		\di = c_v\int\int_{\t_0}^{4^{-1/q}\T}(y^{-q}-4\T^{-q})_+dydx + R\int_{0}^{t}\int\t^{1-q}v^{-1}u_x\left(1-4(\t/\T)^q\right)_+dxds \\[3mm]
 		\di\quad+ 4q\tilde{\kappa}\int_{0}^{t}\int_{\{\t<4^{-1/q}\T\}}v^{-1}\t^\beta\t_x\T^{-q-1}\T_xdxds,
 	\end{array}
 \end{equation}
 which yield that
 \begin{equation}
 	\begin{array}{ll}
 		 \di\int\int_{\t}^{4^{-1/q}\T}(y^{-q}-4\T^{-q})_+dydx + q\tilde{\kappa}\int_{0}^{t}\int_{\t<4^{-1/q}\T}v^{-1}\t^{\beta-q-1}\zeta_x^2dxds\\[3mm]
 		\di\quad +\tilde{\mu} \int_{0}^{t}\int v^{-1}\psi_x^2(\t^{-q}-4\T^{-q})_+dxds+\tilde{\mu} \int_{0}^{t}\int v^{-1}U_x^2(\t^{-q}-4\T^{-q})_+dxds\\[3mm]
 		\di \leq \int\int_{\t_0}^{4^{-1/q}\T}(y^{-q}-4\T^{-q})_+dydx +C\int_{0}^{t}\int_{\t<4^{-1/q}\T}v^{-1}\t^{\beta-q-1}|\zeta_x\T_x|dxds\\[3mm]
 		 \di\quad + C\int_{0}^{t}\int v^{-1}|\psi_x U_x|(\t^{-q}-4\T^{-q})_+dxds +C\int_{0}^{t}\int_{\t<4^{-1/q}\T}v^{-1}\t^\beta\t_x\T^{-q-1}\T_xdxds\\[3mm]
 		 \di\quad + \varepsilon\int_{0}^{t}\int(\psi_x^2+U_x^2)\t^{-q}dx +C(\varepsilon)\int_{0}^{t}\int_{\t<4^{-1/q}\T}\t^{1-q}dx\max_{x\in\mathbb{R}}(1-4(\t/\T)^q)_+^2ds.

\label{beta-1new}
 	\end{array}
 \end{equation}

Next, direct computation yields that
\begin{equation}
	\begin{array}{ll}
		\di\int\int_{\t}^{4^{-1/q}\T}(y^{-q}-4\T^{-q})_+dydx\geq\frac{1}{q-1}\int_{\{\t<4^{-1/q}\T\} }\t^{1-q}dx -C_0.
		 \label{beta-2}
	\end{array}
\end{equation}
And that
\begin{equation}
	\begin{array}{ll}
		\label{beta-3}
\di\int_{0}^{t}\int v^{-1}\psi_x^2(\t^{-q}-4\T^{-q})_+dxds\\[3mm]
\di\quad= \int_{0}^{t}\int v^{-1}\psi_x^2\t^{-q}(1-4(\t/\T)^{q})_+dxds\\[3mm]
\di\quad \geq C^{-1}\int_{0}^{t}\int_{\{\t<1/2\T\}}\t^{-q}\psi_x^2dxds\\[3mm]
\di\quad\geq C^{-1}\int_{0}^{t}\int\t^{-q}\psi_x^2dxds -C\int_{0}^{t}\int_{\{\t>1/4\T\}}\t^{-1}\psi_x^2dxds\\[3mm]
\di\quad\geq C^{-1}\int_{0}^{t}\int\t^{-q}\psi_x^2dxds-C_0.
	\end{array}
\end{equation}

For the right hand of (\ref{beta-1new}), we have
\begin{equation}
	\begin{array}{ll}
		\di\int_{0}^{t}\int_{\{\t<4^{-1/q}\T\}}v^{-1}\t^{\beta-q-1}|\zeta_x\T_x|dxds\\[3mm]
		\di\quad\leq\varepsilon\int_{0}^{t}\int\t^{\beta-q-1}\zeta_x^2dxds + C(\varepsilon)\int_{0}^{t}\int_{\{\t<4^{-1/p}\T\}}\t^{\beta-q-1}\T_x^2dxds\\[3mm]
		\di\quad\leq \varepsilon\int_{0}^{t}\int\t^{\beta-q-1}\zeta_x^2dxds +C(\varepsilon,M)\int_{0}^{t}\int_{\{\t<4^{-1/p}\T\}}\zeta^2\T_x^2dxds\\[3mm]
		\di\quad\leq \varepsilon\int_{0}^{t}\int\t^{\beta-q-1}\zeta_x^2dxds +C_0,
		\label{beta-4}
	\end{array}
\end{equation}
Similarly,
\begin{equation}
	\begin{array}{ll}
		\label{beta-5}
		\di\int_{0}^{t}\int v^{-1}|\psi_x U_x|(\t^{-q}-4\T^{-q})_+dxds\leq\varepsilon\int_{0}^{t}\int \psi_x^2(\t^{-q}-4\T^{-q})_+dxds + C_0,
	\end{array}
\end{equation}
and
\begin{equation}
	\begin{array}{ll}
		\label{beta-6}
		\di\int_{0}^{t}\int_{\{\t<4^{-1/q}\T\}}v^{-1}\t^\beta\t_x\T^{-q-1}\T_xdxds\\[3mm]
		\di\quad\leq C\int_{0}^{t}\int_{\{\t<4^{-1/q}\T\}}\frac{\t^\beta}{\t^2}\t_x^2dxds + C\int_{0}^{t}\int_{\{\t<4^{-1/q}\T\}}\t^{\beta+2}\T^{-2q-2}\T_x^2dxds\\[3mm]
		\di\quad \leq C\int_{0}^{t}\int\frac{\t^\beta}{\t^2}\zeta_x^2dxds+ C\int_{0}^{t}\int_{\{\t<4^{-1/q}\T\}}\left(\t^{\beta+2}\T^{-2q-2}\T_x^2 + \frac{\t^\beta}{\t^2}\T_x^2\right)dxds\\[3mm]
		\di\quad\leq \int_{0}^{t}\int\frac{\t^\beta}{\t^2}\zeta_x^2dxds +C(M)\int_{0}^{t}\int_{\{\t<4^{-1/q}\T\}}\zeta^2\T_x^2dxds\\[3mm]
		\di\quad\leq \int_{0}^{t}\int\frac{\t^\beta}{\t^2}\zeta_x^2dxds +C \leq C_0.

	\end{array}
\end{equation}

Moreover, it follows from (\ref{basic}) and (\ref{LL2}) yield that

\begin{equation}
	\begin{array}{ll}
		\label{beta-7}
		\di\int_{0}^{t}\max_{x\in\mathbb{R}}(1-4(\t/\T)^q)_+^2ds\leq\int_{0}^{t}\left(\max_{x\in\mathbb{R}}\int_{x}^{\infty}\left(-[1-4(\t/\T)^q]\right)_ydy\right)^2ds\\[3mm]
		\di\quad\leq \int_{0}^{t}\left( \int_{\{\t<4^{-1/q}\T\}}q(\t/\T)^{q-1} \frac{\zeta_y\T-\zeta\T_y}{\T^2}dy\right)^2ds\\[3mm]
		\di\quad\leq \int_{0}^{t}\left( 
		C\left(\int_{\{\t<4^{-1/q}\T\}}\t^{-1+\beta/2}\zeta_ydx\right)^2 + C(M) \left(\int_{\{\t<4^{-1/q}\T\}}\T_y\zeta dy\right)^2\right)ds\\[3mm]
		\di\quad\leq C\int_{0}^{t}\int\t^{\beta-2}\zeta_y^2dxds +C(M)\int\T_y^2\zeta^2dxds\leq C_0,
	\end{array}
\end{equation}

and  
\begin{equation}
	\begin{array}{ll}
		\label{beta-8}
		\di\int_{0}^{t}\int\t^{\beta-q-1}\zeta_x^2dxdt\\[3mm]
		\di\quad\leq C\int_{0}^{t}\int_{\{\t<4^{-1/q}\T\}}v^{-1}\t^{\beta-q-1}\zeta_x^2dxds + C\int_{0}^{t}\int_{\{\t>1/2\T\}}\t^{\beta-2}\zeta_x^2dxdst\\[3mm]
		\di\quad\leq C\int_{0}^{t}\int_{\{\t<4^{-1/q}\T\}}v^{-1}\t^{\beta-q-1}\zeta_x^2dxds + C_0.
	\end{array}
\end{equation}
Combining(\ref{beta-1new})-(\ref{beta-8}) and Gronwall's inequality gives 

\begin{equation}
	\di\int_{0}^{t}\int\frac{\t^\beta\zeta_x^2}{\t^{q+1}}dxds + \int_{0}^{t}\int\frac{\psi_x^2}{\t^q}dxds \leq C(q)
			\label{theta^beta}
\end{equation}
for $q\geq\beta+4$. 

Obviously, it follows from (\ref{basic}) that (\ref{theta^beta}) holds for $q=1$. Meanwhile, (\ref{theta^beta}) holds for all $q\geq 1$. We choose $q=\beta$ in (\ref{theta^beta}) gives
\begin{equation}
	\int_{0}^{T}\int\t^{-1}\zeta_x^2dxdt\leq C_0,
\end{equation}
which means that lemma\ref{beta>1} remains to prove
\begin{equation}
	\int_{0}^{T}\int\psi_x^2dxdt\leq C_0.
	\label{psx2}
\end{equation}
Indeed, multiplying(\ref{ener})  by $(\t-2\T)_+\t^{-1}$ and integrating the resulting equality $\mathbb{R} \times (0,t)$ yields
\begin{equation}
	\begin{array}{ll}
		\di \tilde{\mu}\int_{0}^{t}\int\frac{\psi_x^2}{v}(\t-2\T)_+\t^{-1}dxds + \tilde{\mu}\int_{0}^{t}\int\frac{U_x^2}{v}(\t-2\T)_+\t^{-1}dxds\\[3mm]
		\di\quad= -\int_{0}^{t}\int\frac{\psi_xU_x}{v}(\t-2\T)_+\t^{-1}dxds +\tilde{\kappa} \int_{0}^{t}\int_{\Omega_2}\frac{-2\t^{\beta}\t_x\T_x}{v\t} + \frac{2\t^\beta\t_x^2\T}{v\t^2} dxds\\[3mm]
		\di\quad\quad+R\int_{0}^{t}\int\frac{(\t-2\T)_+}{v}u_xdxds + c_v\int\left(\int_{2\T}^{\t(t)}(y-2)_+y^{-1}dy-\int_{2\T}^{\t_0}(y-2)_+y^{-1}dy\right)dx\\[3mm]
		\di\quad\leq \varepsilon\int_{0}^{t}\int\frac{\psi_x^2}{v}(\t-2\T)_+\t^{-1}dxds + C(\varepsilon)\int_{0}^{t}\int\frac{U_x^2}{v}(\t-2\T)_+\t^{-1}dxds \\[3mm]
		\di\quad\quad + C\int_{0}^{t}\int_{\Omega_2}\frac{\t^\beta\zeta_x^2}{v\t^2}\T dxds + C\int_{0}^{t}\int_{\Omega_2}\frac{\t^\beta\T_x^2}{v\t^2}\T dxds + C\int_{0}^{t}\int_{\Omega_2}\left|\frac{\t^\beta\t_x\T_x}{v\t}\right|dxds\\[3mm]
		\di\quad\quad +C \int_{0}^{t}\int\frac{(\t-2\T)_+}{v}\psi_x dxds + C\int_{0}^{t}\int\frac{(\t-2\T)_+}{v}U_xdxds  +C_0 \\[3mm]
		\di\quad\leq \varepsilon\int_{0}^{t}\int\frac{\psi_x^2}{v}(\t-2\T)_+\t^{-1}dxds + C(M,\varepsilon)\int_{0}^{t}\int_{\Omega_2}\zeta^2(U_x^2+\T_x^2)dxds \\[3mm]
		\di\quad\quad+\varepsilon\int_{0}^{t}\int\psi_x^2dxds + C(\varepsilon)\int_{0}^{t}\max_{x\in\mathbb{R}}(\t-3/2\T)_+^{\beta+1}dt + C_0
		\\[3mm]
		\di\quad\leq \varepsilon\int_{0}^{t}\int\frac{\psi_x^2}{v}(\t-2\T)_+\t^{-1}dxds + \varepsilon\int_{0}^{t}\int\psi_x^2dxds + C(\varepsilon)\int_{0}^{t}\max_{x\in\mathbb{R}}(\t-3/2\T)_+^{\beta+1}dt + C_0.
		\label{psx21}
	\end{array}
\end{equation}

Then, direct calculation shows that for $\chi\geq-1$,
\begin{equation*}
	\begin{array}{ll}
		\di\max_{x\in\mathbb{R}}(\t-3/2\T)_+^{\chi+3} = \max_{x\in\mathbb{R}}\int_{x}^{\infty}\left(-\partial_x(\t-3/2\T)_+^{\chi+3}\right)dx\\[3mm]
		\di\quad\leq C \int_{\{\t>3/2\T\}}(\t-3/2\T)_+^{\chi+2}|\t_x-3/2\T_x|dx\\[3mm]
		\di\quad\leq \left(C\left(\int_{\{\t>3/2\T\}}\t_x^2\t^\chi dx\right)^{1/2} + C\left(\int_{\{\t>3/2\T\}}\T_x^2\t^\chi dx\right)^{1/2}\right)\max_{x\in\mathbb{R}}(\t-3/2\T)_+^{(\chi+3)/2},
	\end{array}
\end{equation*} 
where in the last inequality we have used (\ref{cedu}). This gives 
\begin{equation}
	\begin{array}{ll}
		\label{psx22}
		\di\max_{x\in\mathbb{R}}(\t-3/2\T)_+^{\chi+3} \leq C \int_{\{\t>3/2\T\}}\zeta_x^2\t^\chi + C\int_{\{\t>3/2\T\}}\T_x^2\t^\chi dx.
	\end{array}
\end{equation}
In particular, since $\beta>1$. choosing $\chi=\beta-2$ in (\ref{psx22}) gives  
\begin{equation}
	\label{psx23}
	\int_{0}^{t}\max_{x\in\mathbb{R}}(\t-3/2\T)_+^{\beta+1}ds\leq C\int_{0}^{t}\int_{\{\t>3/2\T\}}\zeta_x^2\t^{\beta-2}dxds + C(M)\int_{0}^{t}\int_{\{\t>3/2\T\}}\zeta^2\T_x^2dxds \leq C_0.
\end{equation} 

Finally, it follows from (\ref{basic}), (\ref{psx21}) and (\ref{psx23}) that
\begin{equation}
	\begin{array}{ll}
		\di \int_{0}^{t}\int\psi_x^2dxds &\di \leq C\int_{0}^{t}\int_{\Omega_3}\frac{\psi_x^2}{v}(\t-2\T)_+\t^{-1}dxds\\[3mm]
		&\di\quad + C\int_{0}^{t}\int_{\t<4\T}\frac{\psi_x^2}{v\t}dxds\\[3mm]
		&\di \leq C_0,
		\label{psx24}
	\end{array}
\end{equation}
which gives (\ref{psx2}) and finish the proof of this lemma.

\

\begin{lemma}\label{high-order}
Suppose that $(\phi,\psi,\zeta)\in X([0,T])$ satisfies $\d=|\t_+-\t_-|\leq\d_0$
with suitable small $\d_0$, it holds
\begin{equation}\label{high-deri}
\di\sup_{0\leq t\leq T}\int(\phi_x^2+\psi_x^2)dx
+\int_0^T\int(\zeta_x^2+(1+\t)\phi_x^2+\psi_{xx}^2)dxdt\leq C_0.
\end{equation}
\end{lemma}
\textbf{Proof}:\quad Due to \eqref{basic}, \eqref{v},  \eqref{new-nergy} and (\ref{beta-new}), some terms of \eqref{phi-x1} can be considered more carefully, that is,
\begin{equation}
\begin{array}{ll}
\di \left|\frac{R}{v}\zeta_x\frac{\tilde v_x}{\tilde v}\right|
+\left|\psi_xU_x\left(\frac{1}{v}-\frac{1}{V}\right)\right|+\frac{\psi_x^2}{v}\leq \frac{R\t}{4v}\left(\frac{\tilde v_x}{\tilde v}\right)^2+
C_0\frac{\zeta_x^2}{\t}+C_0\psi_x^2+C_0\phi^2U_x^2.\\[3mm]
\end{array}
\end{equation}
The other terms in \eqref{phi-x1} can be estimated the same as in step 2 in Lemma \ref{basic-lemma}.
Integrating \eqref{phi-x1} over $\mathbb{R}\times(0,t)$, recalling \eqref{basic} and \eqref{new-nergy}, we have
\begin{equation}
	\label{phix-p}
	\begin{array}{ll}
&\di \sup_{0\leq t\leq T}\int\phi_x^2 dx+C\int_0^T\int(1+\t)\phi_x^2dxdt\\[3mm]
&\di\quad\leq \sup_{0\leq t\leq T}\int\phi_x^2 dx+\frac{1}{2}\int_0^T\int(\T+\t)\phi_x^2dxdt\\[3mm]
&\di\quad\leq C_0 + \frac{1}{2}\int_0^T\int(\T-\t)_+\phi_x^2dxdt\\[3mm]
&\di\quad\leq C_0+\frac{1}{4}\int_{0}^{T}\int\phi_x^2dxdt + C\int_{0}^{T}\max_{x\in\mathbb{R}}(\T-\t)_+^4\int\phi_x^2dxdt.		
	\end{array}
\end{equation}
Then, it follows from (\ref{new-nergy}), (\ref{beta-new}) and (\ref{cedu}) that
\begin{equation}
	\begin{array}{ll}
		\di\int_{0}^{t}\max_{x\in\mathbb{R}}(\T-\t)_+^4 ds \leq C\int_{0}^{t}\left(\max_{x\in\mathbb{R}}\int_{x}^{\infty}(1-\frac{\t}{\T})_+\left(\frac{\t}{\T}\right)_xdx\right)^2ds\\[4mm]
		\di\quad\leq C\int_{0}^{t}\left(\int(\T-\t)_+\T^{-1}\left(\frac{\T\zeta_x-\T_x\zeta}{\T^2}\right)dx\right)^2ds\\[4mm]
		\di\quad\leq C\int_{0}^{t}\left(\int(\T-\t)_+\t^{-1/2}|\zeta_x|dx\right)^2ds + C\int_{0}^{t}\left(\int(\T-\t)_+\T^{-3}|\T_x\zeta| dx\right)^2ds\\[3mm]
		\di\quad\leq C\int_{0}^{t}\left(\int(\T-\t)_+^2dx\int\t^{-1}\zeta_x^2dx\right)ds + C(M)\int_{0}^{t}\int\T_x^2\zeta^2dxds \leq C_0,
		
	\end{array}
\end{equation}
which together with (\ref{phix-p}) and Gronwall's inequality gives
\begin{equation}
	\begin{array}{ll}
		\di \sup_{0\leq t\leq T}\int\phi_x^2 dx+C\int_0^T\int(1+\t)\phi_x^2dxdt\leq C_0.
		\label{phix}
	\end{array}
\end{equation}

Multiplying $\eqref{perturb}_2$ by $-\psi_{xx}$, integrating the resulted equation over $\mathbb{R}\times(0,t)$,
and noticing that
\begin{equation*}
\begin{array}{ll}
\di(p-p_+)_x&\di=\left(\frac{R\zeta-p_+\phi}{v}\right)_x=\frac{R\zeta_x-p_+\phi_x}{v}
-\frac{R\zeta-p_+\phi}{v^2}\phi_x-\frac{R\zeta-p_+\phi}{v^2}V_x\\[3mm]
&\di=\frac{R\zeta_x}{v}-\frac{R\t\phi_x}{v^2}-\frac{R\zeta-p_+\phi}{v^2}V_x.
\end{array}
\end{equation*}
Then we have
\begin{equation}
\begin{array}{ll}\label{6.1}
\di\int\frac{\psi_x^2}{2}dx+\tilde{\mu}\int_0^t\int\frac{\psi_{xx}^2}{v}dxds
=\int\frac{\psi_{0x}^2}{2}dx+\int_0^t\int\left(\frac{R\zeta_x}{v}-\frac{R\t\phi_x}{v^2}
-\frac{R\zeta-p_+\phi}{v^2}V_x\right)\psi_{xx}dxds\\
\di\quad-\tilde{\mu}\int_0^t\int\psi_x\left(\frac{1}{v}\right)_x\psi_{xx}dxds
+\tilde{\mu}\int_0^t\int\left(\frac{U_x}{V}-\frac{U_x}{v}\right)_x\psi_{xx} dxds
+\int_0^t\int\widetilde R_1\psi_{xx}dxds.
\end{array}
\end{equation}

Next, we will estimate each term on the right-hand side of \eqref{6.1}.
From \eqref{important}, \eqref{new-nergy}, \eqref{beta-new} and \eqref{phix}, one has
\begin{equation}
\begin{array}{ll}\label{psi-xx-1}
\di\left|\int_0^t\int\left(\frac{R\zeta_x}{v}-\frac{R\t\phi_x}{v^2}
-\frac{R\zeta-p_+\phi}{v^2}V_x\right)\psi_{xx}dxds\right|\\
\di\leq \frac{\tilde{\mu}}{16}\int_0^t\int\frac{\psi_{xx}^2}{v}dxds
+C_0\int_0^t\int\zeta_x^2+\theta^2\phi_x^2+(\phi^2+\zeta^2)V_x^2dxds\\
\di\leq \frac{\tilde{\mu}}{16}\int_0^t\int\frac{\psi_{xx}^2}{v}dxds+C_0\int_0^t\int\zeta_x^2dxds\\[3mm]
\di\quad +C_0\int_0^t\left(\max_{x\in \mathbb{R}}\left(\t-\frac{3}{2}\T\right)_+^2+1\right)\int\phi_x^2dxds+C_0\\
\di\leq \frac{\tilde{\mu}}{16}\int_0^t\int\frac{\psi_{xx}^2}{v}dxds+C_0\int_0^t\int\zeta_x^2dxds
+C_0\int_{0}^{t}\max_{x\in \mathbb{R}}\left(\t-\frac{3}{2}\T\right)_+^2ds +C_0.
\end{array}
\end{equation}
Noticing that (\ref{psx22}) where we choose $\chi=-1$, combining (\ref{new-nergy})  and (\ref{beta-new}) gives 
\begin{equation}
	\di\int_{0}^{T}\max_{x\in \mathbb{R}}(\t-\frac{3}{2}\T)_+^2dt \leq C_0.
	\label{delta=-1}
\end{equation}

Next, recalling \eqref{basic}, \eqref{new-nergy}, \eqref{beta-new} and 
\eqref{phix},  we obtain by Cauchy's inequality and Sobolev's inequality that 
\begin{equation}
\begin{array}{ll}\label{psi-xx-2}
\di\left|\tilde{\mu}\int_0^t\int\psi_x\left(\frac{1}{v}\right)_x\psi_{xx}dxds\right|
\leq C\int_0^t\int\left|\frac{\psi_x\phi_x\psi_{xx}}{v^2}\right|+\left|\frac{\psi_xV_x\psi_{xx}}{v^2}\right|dxds\\
\di\leq\frac{\tilde{\mu}}{16}\int_0^t\int\frac{\psi_{xx}^2}{v}dxds
+C_0\int_0^t\int \psi_x^2\phi_x^2+\psi_x^2\Theta_x^2dxds\\
\di\leq\frac{\tilde{\mu}}{16}\int_0^t\int\frac{\psi_{xx}^2}{v}dxds
+C_0\int_0^t\|\psi_x\|^2_{L^{\infty}}\|\phi_x\|^2ds
+C_0\int_0^t\int\psi_x^2dxds\\
\di\leq\frac{\tilde{\mu}}{16}\int_0^t\int\frac{\psi_{xx}^2}{v}dxds
+C_0\int_0^t\|\psi_x\|\|\psi_{xx}\|ds
+C_0\int_0^t\int\psi_x^2dxds\\
\di\leq\frac{\tilde{\mu}}{8}\int_0^t\int\frac{\psi_{xx}^2}{v}dxds
+C_0\int_0^t\int\psi_x^2dxds\\
\di\leq\frac{\tilde{\mu}}{8}\int_0^t\int\frac{\psi_{xx}^2}{v}dxds+C_0.
\end{array}
\end{equation}
Similarly,
\begin{equation}
\begin{array}{ll}\label{psi-xx-3}
\di\tilde{\mu}\int_0^t\int\left(\frac{U_x}{V}-\frac{U_x}{v}\right)_x\psi_{xx} dxds\\
\di=\tilde{\mu}\int_0^t\int\left(\frac{U_{xx}}{V}-\frac{U_{xx}}{v}-\frac{U_xV_x}{V^2}+\frac{v_xU_x}{v^2}\right)\psi_{xx}dxds\\
\di=\tilde{\mu}\int_0^t\int\left(\frac{\phi U_{xx}}{vV}-U_xV_x\frac{\phi(\phi+2V)}{v^2V^2}+\frac{\phi_xU_x}{v^2}\right)\psi_{xx}dxds\\
\di\leq\frac{\tilde{\mu}}{16}\int_0^t\int\frac{\psi_{xx}^2}{v}dxds
+C_0\int_0^t\int(\phi^2(U_{xx}^2+\Theta_x^2U_x^2)+\phi^4U_x^2\Theta_x^2+\phi_x^2U_x^2)dxds\\
\di\leq\frac{\tilde{\mu}}{16}\int_0^t\int\frac{\psi_{xx}^2}{v}dxds
+C_0+C(C_0,M)\d^2\int_0^t\int\phi_x^2dxds\\
\di\leq\frac{\tilde{\mu}}{16}\int_0^t\int\frac{\psi_{xx}^2}{v}dxds+C_0,
\end{array}
\end{equation}
and
\begin{equation}
\begin{array}{ll}\label{psi-xx-4}
\di\left|\int_0^t\int\widetilde R_1\psi_{xx}dxds\right|
\leq\frac{\tilde{\mu}}{16}\int_0^t\int\frac{\psi_{xx}^2}{v}dxds+C_0\int_0^t\int\widetilde R_1^2dxds\\
\di\leq\frac{\tilde{\mu}}{16}\int_0^t\int\frac{\psi_{xx}^2}{v}dxds+C_0\d^2\int_0^t\int(1+s)^{-3}e^{-\frac{c_1x^2}{1+s}}dxds\\
\di\leq\frac{\tilde{\mu}}{16}\int_0^t\int\frac{\psi_{xx}^2}{v}dxds+C_0.
\end{array}
\end{equation}
Substituting \eqref{psi-xx-1}-\eqref{psi-xx-4} into \eqref{6.1} shows
\begin{equation}\label{psi-xx}
\di\sup_{0\leq t\leq T}\int\psi_x^2dx+\int_0^T\int\psi_{xx}^2dxdt
\leq C_0+C_0\int_{0}^{t}\int\zeta_x^2dxdt.
\end{equation}

Finally, we estimate the right hand of (\ref{psi-xx}). On the one hand, if $\beta\geq 2$, choosing $q=\beta-1$ in (\ref{theta^beta}) gives
\begin{equation}
	\begin{array}{ll}
		\di\int_{0}^{T}\int\zeta_x^2dxdt\leq C_0	,\end{array}
\end{equation}
which along with (\ref{psi-xx}) shows
\begin{equation}
	\di\sup_{0\leq t\leq T}\int\psi_x^2dx+\int_0^T\int\psi_{xx}^2dxdt
	+\int_{0}^{T}\int\zeta_x^2dxdt\leq C_0.
	\label{zetax-f}
\end{equation}
On the other hand, if $\beta\in (0,2)$, multiplying (\ref{ener}) by $(\t-2\T)_+\t^{-\beta/2}$ and integration by parts gives 
\begin{equation}
	\begin{array}{ll}
		\di\tilde{\kappa}\int_{0}^{t}\int\frac{\t^{\beta/2}\zeta_x^2(2-\beta)}{2v}dxds +\tilde{\kappa} \int_{0}^{t}\int\frac{\beta\t^{-1+\beta/2}\T\zeta_x^2}{v}dxds\\[3mm]
		\di \quad\leq \tilde{\kappa}\int_{0}^{t}\int \frac{2\t^{\beta/2}\t_x\T_x}{v} - c_v\int\int_{2\T}^{\t(t)}(y-2\T)_+y^{-\beta/2}dydx + c_v\int\int_{2\T}^{\t(t)}(y-2\T)_+y^{-\beta/2}dydx\\[4mm]
		\di \quad \quad -R\int_{0}^{t}\int\frac{(\t-2\T)_+\t^{1-\beta/2}u_x}{v}dxds + \tilde{\mu}\int_{0}^{t}\frac{(\t-2\T)_+\t^{-\beta/2}u_x^2}{v}dxds\\[4mm]
		\di \quad \quad - \tilde{\kappa}\int_{0}^{t}\int_{\Omega_2}\frac{\t^{\beta/2}|\zeta_x||\T_x|(2-\beta)}{v}dxds -\tilde{\kappa} \int_{0}^{t}\int_{\Omega_2}\frac{\beta\t^{-1+\beta/2}\T|\zeta_x||\T_x|}{v}dxds\\[3mm]
		\di\quad\leq C_0+ C\int_{0}^{t}\frac{\t^\beta\zeta_x^2}{\t^2v}dxds + C(M)\int_{0}^{t}\int_{\Omega_2}\T_x^2\zeta^2dxds + C\int_{0}^{t}\int(\t-2\T)_+\t^{2-\beta/2}dxds \\[3mm]
		\di\quad\quad + 2\int_{0}^{t}\int(\t-2\T)_+\t^{-\beta/2}u_x^2dxds\\[3mm]
		\di \quad\leq C_0 + C\int_{0}^{t}\left(\max_{x\in \mathbb{R}}(\t-\frac{3}{2}\T)_+^2\int_{\Omega_2}\t^{1-\beta/2}dx\right)ds \\[3mm]
		\di\quad\quad + 4\int_{0}^{t}\int(\t-2\T)_+\psi_x^2dxds + 4\int_{0}^{t}\int_{\Omega_2}(\t-2\T)_+U_x^2dx\\[3mm]
		\di \quad\leq C_0 + C\int_{0}^{t}\max_{x\in \mathbb{R}}(\t-\frac{3}{2}\T)_+^2ds + C\int_{0}^{t}(\int\psi_x^2dx)^2ds + C(M)\int_{0}^{t}\int_{\Omega_2}\zeta^2U_x^2dxds,
	\end{array}
	\label{psil2l4}
\end{equation}
where in the last inequality we have used (\ref{omega-bound}) and (\ref{cedu}). Since $\beta<2$, it follows from (\ref{psil2l4}) ,(\ref{delta=-1}) and (\ref{basic}) that
\begin{equation}
	\begin{array}{ll}
		\di\int_{0}^{t}\int\zeta_x^2dxds &\di\leq C\int_{0}^{t}\int_{\zeta<2\T}\t^{\beta-2}\zeta_x^2 + C\int_{0}^{t}\int_{\Omega_2}\t^{\beta/2}\zeta_x^2 \\[3mm]
		&\di\leq C_0 + C\int_{0}^{t}\left(\int\psi_x^2dx\right)^2dt,
		\label{beta0-2}
	\end{array}
\end{equation}
which with (\ref{psi-xx}), (\ref{new-nergy}), (\ref{beta-new}) and Gronwall's inequality shows that (\ref{zetax-f}) still holds.

Combining (\ref{phix-p}), (\ref{phix}), (\ref{zetax-f}) and  (\ref{beta0-2}) finishes lemma\ref{high-order}.

%
\
%
%
%
%

\begin{lemma}
	\label{theta-supbel}
Suppose that $(\phi,\psi,\zeta)\in X([0,T])$ satisfies $\d=|\t_+-\t_-|\leq\d_0$
with suitable small $\d_0$,	there exists a positive constant $C_0$ such that for any $(x,t)\in \mathbb{R} \times [0,T]$
	\begin{equation}
		\di C_0^{-1}\leq \t(x,t)\leq C_0.
		\label{theta-fina}
	\end{equation}
	
\end{lemma}
\textbf{Proof}:\quad For $q>\beta+1$, multiplying ($\ref{ener}$) by $(\t-2\T)_+^{q-1}$ and integrating the resultant equality in $\mathbb{R}\times(0,t)$ leads to 
\begin{equation}
	\begin{array}{ll}
		\di(q-1)\tilde{\kappa}\int_{0}^{t}\int\frac{\t^\beta(\t-2\T)_+^{q-2}\zeta_x^2}{v}dxds\\[3mm]
		\di\quad\leq C\int_{0}^{t}\int\frac{\t^\beta(\t-2\T)_+^{q-2}\T_x^2}{v}dxds + C\int_{0}^{t}\int\frac{\t^\beta(\t-2\T)_+^{q-2}|\t_x||\T_x|}{v}dxds \\[4mm]
		\di\quad\quad+ c_v\int\int_{2\T}^{\t(t)}(y-2\T)_+^{q-1}dydx +c_v \int\int_{2\T}^{\t_0}(y-2\T)_+^{q-1}dydx\\[4mm]
		\di\quad\quad -\tilde{\mu}\int_{0}^{t}\int\frac{\t(\t-2\T)_+^{q-1}u_x}{v}dxds + \tilde{\mu}\int_{0}^{t}\int\frac{(\t-2\T)_+^{q-1}u_x^2}{v}dxds\\[4mm]
		\di\quad\leq
		C\varepsilon\int_{0}^{t}\int\frac{\t^\beta(\t-2\T)_+^{q-2}\zeta_x^2}{v}dxds +  C(M,\varepsilon)\int_{0}^{t}\int_{\zeta>\T}\zeta^2(\T_x^2+U_x^2)dxds + C_0\\[4mm]
		\di\quad\quad   +C(\varepsilon)\int_{0}^{t}\int\frac{(\t-2\T)_+^{q-1}\psi_x^2}{v}dxds + \varepsilon\int_{0}^{t}\int\frac{\t^2(\t-2\T)_+^{q-1}}{v}dxds\\[4mm]
		\di\quad\leq
		C\varepsilon\int_{0}^{t}\int\frac{\t^\beta(\t-2\T)_+^{q-2}\zeta_x^2}{v}dxds + C_0 +C(\varepsilon)\int_{0}^{t}\left(\int\psi_x^2dx\right)^{(\beta+q+1)/(\beta+2)}ds \\[4mm]
		\di\quad\quad + C\varepsilon\int_{0}^{t}\max_{x\in \mathbb{R}}(\t-\frac{3}{2}\T)_+^{\beta+q+1}ds\\[4mm]
		\di\quad\leq C\varepsilon\int_{0}^{t}\int\frac{\t^\beta(\t-2\T)_+^{q-2}\zeta_x^2}{v}dxds + C_0 +  C(\varepsilon)\int_{0}^{t}\left(\int\psi_x^2dx\right)^{(\beta+q+1)/(\beta+2)}ds\\[4mm]
		\di\quad\quad  + C\varepsilon\int_{0}^{t}\int_{\t>3/2\T}\zeta_x^2\t^{\beta+q-2}dxds  + C(M)\int_{0}^{t}\int_{\t>3/2\T}\zeta^2\T_x^2dxds,
			\label{int1}
	\end{array}
\end{equation}
where in the last equality we have used (\ref{zetax-f}),  (\ref{psx22}) and lemma2.1.

Then,it follows from (\ref{basic}) that
\begin{equation*}
	\begin{array}{ll}
	\di\int_{0}^{t}\int\t^{\beta+q-2}\zeta_x^2dxds&\di\leq C\int_{0}^{t}\int\frac{\t^\beta(\t-2\T)_+^{q-2}\zeta_x^2}{v}dxds + C\int_{0}^{t}\int_{\t<3\T}\t^{\beta-2}\zeta_x^2dxds\\[3mm]
	&\di\leq C\int_{0}^{t}\int\frac{\t^\beta(\t-2\T)_+^{q-2}\zeta_x^2}{v}dxds+C_0,
	\end{array}
\end{equation*}
which together with (\ref{int1}) and (\ref{new-nergy}),  (\ref{beta-new}) and (\ref{high-deri})gives
\begin{equation}
	\begin{array}{ll}
		\di\int_{0}^{T}\int\t^{\beta+q-2}\zeta_x^2dxdt \leq C_0.
		\label{int-t+}
	\end{array}
\end{equation}

Multiplying $\eqref{perturb}_3$ by $\theta^\beta\zeta_t$, then integrating the resulted equation over
$\mathbb{R}$, we have
\begin{equation}
	\begin{array}{ll}
		\di c_v\int\t^\beta\zeta_t^2dx &\di= \tilde{\kappa}\int\left(\frac{\t^\beta\t_x}{v}-\frac{\T^\beta\T_x}{V}\right)_x\t^\beta\zeta_tdx +  \tilde{\mu}\int\left(\frac{u_x^2}{v}-\frac{U_x^2}{V}\right)\t^\beta\zeta_tdx\\[3mm]
		&\di\quad -\int(\frac{R\t}{v}u_x-\frac{R\T}{V}U_x)\t^\beta\zeta_tdx - \int\widetilde R_2\t^\beta\zeta_tdx\\[3mm]
		&\di= \tilde{\kappa}\int\left(\frac{\t^\beta\zeta_x}{v}\right)_x\t^\beta\zeta_tdx + \tilde{\kappa}\int\left(\frac{\t^\beta\T_x}{v}-\frac{\T^\beta\T_x}{V}\right)_x\t^\beta\zeta_t\\[3mm]
		&\di\quad   +  \tilde{\mu}\int\left(\frac{u_x^2}{v}-\frac{U_x^2}{V}\right)\t^\beta\zeta_tdx -\int(\frac{R\t}{v}u_x-\frac{R\T}{V}U_x)\t^\beta\zeta_tdx- \int\widetilde R_2\t^\beta\zeta_tdx\\[3mm]
		& \di\triangleq \sum_{K=1}^5K_i.
		\label{supthet}
	\end{array}
\end{equation}
We will estimate (\ref{supthet}) one by one. First, we have
\begin{equation}
	\begin{array}{ll}
	\di	K_1 &\di= -\tilde{\kappa}\int\left(\frac{\t^\beta\zeta_x}{v}\right)\left((\t^\beta\zeta_x)_t + (\t^\beta\T_x)_t-(\t^\beta\T_t)_x\right)dx\\[3mm]
	&\di =-\tilde{\kappa}\int\left(\frac{\t^\beta\zeta_x}{v}\right)\left((\t^\beta\zeta_x)_t + \beta\theta^{\beta-1}(\zeta_t\T_x-\zeta_x\T_t)\right)dx\\[3mm]
	&\di = -\tilde{\kappa}\int\frac{1}{2v}\frac{d}{dt}\left(\t^\beta\zeta_x\right)^2dx   -\tilde{\kappa}\int\left(\frac{\t^\beta\zeta_x}{v}\right)\left(\beta\theta^{\beta-1}(\zeta_t\T_x-\zeta_x\T_t)\right)dx\\[3mm]
	&\di = -\tilde{\kappa}\int\left(\frac{\t^{2\beta}\zeta_x^2}{2v}\right)_t + \frac{(\t^\beta \zeta_x)^2v_t}{2v^2}dx-\tilde{\kappa}\int\left(\frac{\t^\beta\zeta_x}{v}\right)\left(\beta\theta^{\beta-1}(\zeta_t\T_x-\zeta_x\T_t)\right)dx.
	\label{k11}
	\end{array}
\end{equation}
For the last two term on the right hand side of (\ref{k11}), it follows from lemma 2.1 that
\begin{equation}
	\begin{array}{ll}
	\di\int\frac{(\t^\beta \zeta_x)^2v_t}{2v^2}dx&\di\leq C\int\frac{\t^{2\beta}\zeta_x^2}{v^2}(\psi_x + U_x)dx	\\[3mm]
	&\di\leq C\max_{x\in \mathbb{R}}\psi_x\int\frac{\t^{2\beta}\zeta_x^2}{v^2}dx + C\int\frac{\t^{2\beta}\zeta_x^2}{v^2}dx\\[3mm]
	&\di\leq C\max_{x\in\mathbb{R}}\psi_x^2 +C \left(\int\frac{\t^{2\beta}\zeta_x^2}{v^2}dx\right)^2 + C\int\frac{\t^{2\beta}\zeta_x^2}{v^2}dx.
	\label{k12}
	\end{array}
\end{equation}
And using (\ref{Theta t}) gives
\begin{equation}
	\begin{array}{ll}
		\di\int\left(\frac{\t^\beta\zeta_x}{v}\right)\left(\beta\theta^{\beta-1}(\zeta_t\T_x-\zeta_x\T_t)\right)dx\\[3mm]
		 \di\quad\leq\frac{c_v}{8}\int\t^\beta\zeta_t^2dx + C\int\t^{3\beta-2}\zeta_x^2\T_x^2dx + \int\frac{\beta\t^{2\beta-1}}{v}\zeta_x^2|\T_t|dx\\[3mm]
		 \di\quad \leq\frac{c_v}{8}\int\t^\beta\zeta_t^2dx + C(M)\delta^2\int\zeta_x^2dx + C\int\frac{\t^{2\beta-1}\zeta_x^2}{v}\left(\left|\frac{\T_{xx}}{V}\right| + \left|\frac{\T_xV_x}{V^2}\right| + \left|\frac{\T U_x}{V}\right| \right)dx\\[3mm]
		 \di\quad\leq \frac{c_v}{8}\int\t^\beta\zeta_t^2dx  + C(M)\delta\int\zeta_x^2dx. 
		 \label{k12}
	\end{array}
\end{equation}

By using (\ref{cedu}) and lemma 2.1, it holds that
\begin{equation}
	\begin{array}{ll}
		\di K_2 &\di = \int\left(\left(\frac{\T_x}{vV}\right)_x(V\t^\beta-v\T^\beta)-\left(\frac{\T_x}{vV}\right)(V\t^\beta-v\T^\beta)_x\right)\t^\beta\zeta_tdx\\[4mm]
		&\di\leq \frac{c_v}{8}\int\t^\beta\zeta_t^2dx + C\int\t^{3\beta}(\T_{xx}^2+\T_x^2v_x^2+\T_x^2V_x^2)dx\\[4mm]
		&\di \quad +C\int\t^\beta\T^{2\beta}(\T_{xx}^2 + \T_x^2v_x^2+\T_x^2V_x^2)dx \\[4mm]
		&\di \quad+C\int\t^{2\beta}\T_x^2V_x^2+\t^{3\beta-2}\T_x^2\t_x^2+\t^{3\beta}\T_x^2v_x^2+\T^{3\beta-2}\T_x^4dx\\[4mm]
		&\di \leq \frac{c_v}{8}\int\t^\beta\zeta_t^2dx+  C(M)\int\left(\T_x^2(\phi_x^2+\zeta_x^2) + \T_{xx}^2 + \T_x^4 + V_x^4\right)dx .
		\label{k2}
	\end{array}
\end{equation} 
Similarly, one has
\begin{equation}
	\begin{array}{ll}
	\label{k3}
	\di K_3&\di=\tilde{\mu}\int\left(\frac{u_x^2}{v}-\frac{U_x^2}{V}\right)\t^\beta\zeta_tdx=\tilde{\mu}\int\left(\frac{\psi_x^2+2\psi_xU_x}{v}-\frac{\phi U_x^2}{vV}\right)\t^\beta\zeta_tdx
	\\[3mm]
	&\di\leq\frac{c_v}{8}\int\t^\beta\zeta_t^2dx + C\int\t^\beta\psi_x^4 + \t^\beta\psi_x^2U_x^2+\t^\beta\phi^2 U_x^4dx\\[3mm]
	&\di\leq\frac{c_v}{8}\int\t^\beta\zeta_t^2dx + C\int\left(\t^{-1}\zeta_x^2+\t^{2\beta-1}\zeta_x^2\right)dx + C\left(\|\psi_x\|^2_{L^2}+\|\psi_{xx}\|^2_{L^2}\right)\\[3mm]
	&\di \quad +C(M)\int(\zeta_x^2\T_x^2+\psi_x^2U_x^2 + \phi^2U_x^4)dx,
	\end{array}
\end{equation}
where in th last equality we used that
\begin{equation*}
	\begin{array}{ll}
		\di\int\t^{\beta}\psi_x^4dx &\di\leq C\int(\t^\beta-2\T^\beta)_+\psi_x^4dx + \int2\T^\beta\psi_x^4dx\\[3mm]
		&\di\leq C\max_{x\in \mathbb{R}}(\t^\beta-2\T^\beta)_+\int\psi_x^4dx + C\int\psi_x^4dx\\[3mm]
		&\di\leq C\int_{-\infty}^{\infty}\left|\beta\t^{\beta-1}\t_x-2\beta\T^{\beta-1}\T_x\right|1_{\{\t>2^{1/\beta}\T\}}dx\int\psi_x^4dx + C\max_{x\in \mathbb{R}}\psi_x^2\int\psi_x^2dx\\[3mm]
		&\di\leq C\left(\int_{\{\t>2^{1/\beta}\T\}}\t^{2\beta-2}\zeta_x^2dx\right)^{1/2}\max_{x\in \mathbb{R}}\psi_x^2\int\psi_x^2dx\\[3mm]
		&\di\quad +\left(C(M)\int_{\{\t>2^{1/\beta}\T\}}\zeta^2\T_x^2dx\right)^{1/2}\max_{x\in \mathbb{R}}\psi_x^2\int\psi_x^2dx ++ C\max_{x\in \mathbb{R}}\psi_x^2\int\psi_x^2dx\\[3mm]
		&\di\leq C\|\psi_x\|^2_{L^2}\|\psi_{xx}\|^2_{L^2}+\int_{\{\t>2^{1/\beta}\T\}}\left(\t^{-1}\zeta_x^2+\t^{2\beta-1}\zeta_x^2\right)dx \\[3mm]
		&\quad\di+C(M)\int_{\{\t>2^{1/\beta}\T\}}\zeta^2\T_x^2dx + C\|\psi_x\|_{L^2}\|\psi_{xx}\|_{L^2} \\[3mm]
		&\di\leq C\int\left(\t^{-1}\zeta_x^2+\t^{2\beta-1}\zeta_x^2\right)dx +C(M)\int\zeta_x^2\T_x^2dx+ C\left(\|\psi_x\|^2_{L^2}+\|\psi_{xx}\|^2_{L^2}\right) .
		
	\end{array}
\end{equation*}

Similarly,
\begin{equation}
	\begin{array}{ll}
		\label{k4}
		\di |K_4| &\di\leq \frac{c_v}{8}\int\t^\beta\zeta_t^2dx + C\int\t^{\beta+2}\psi_x^2 + \zeta^2U_x^2 + \T^2\phi^2U_x^2dx\\[3mm]
		&\di\leq \frac{c_v}{8}\int\t^\beta\zeta_t^2dx +C\int\t^{2\beta+2}\zeta_x^2dx\int\psi_x^2dx + C\int\psi_x^2dx\\[3mm]
		&\di \quad+ C(M)\int(\zeta^2+\phi^2)U_x^2 + \delta\psi_x^2dx,
	\end{array}
\end{equation}
where we used that 
\begin{equation*}
	\begin{array}{ll}
		&\di\int\t^{\beta+2}\psi_x^2dx\di\leq\int\left(\t^{\beta+2}-2\T^{\beta+2}\right)_+\psi_x^2dx + \int2\T^{\beta+2}\psi_x^2dx\\[3mm]
		&\quad\di\leq\int_{-\infty}^{+\infty}|(\beta+2)\t^{\beta+1}\t_x-2(\beta+2)\T^{\beta+1}\T_x|1_{\{\t>2^{1/(\beta+2)}\T\}}dx\int\psi_x^2dx + C\int\psi_x^2dx\\[3mm]
		&\quad\di\leq C\left(\int_{-\infty}^{+\infty}\t^{\beta+1}\zeta_x1_{\{\t>2^{1/(\beta+2)}\T\}}dx + C(M)\int|\T_x|dx\right)\int\psi_x^2dx+C\int\psi_x^2dx\\[3mm]
		&\quad\di\leq C\left(\left(\int_{-\infty}^{+\infty}\t^{2\beta+2}\zeta_x^2dx\right)^{1/2} + \delta C(M)dx\right)\int\psi_x^2dx + C\int\psi_x^2dx\\[3mm]
		&\quad\di\leq C\int\t^{2\beta+2}\zeta_x^2dx\int\psi_x^2dx +C(M)\delta\int\psi_x^2dx + C\int\psi_x^2dx.
	\end{array}
\end{equation*}

Recalling lemma 2.1, one has
\begin{equation}
	\begin{array}{ll}
		\label{k5}
		\di \int_{0}^{t}|K_5|ds&\di\leq\frac{c_v}{8}\int_{0}^{t}\int\t^\beta\zeta_t^2dxdx +C_0\delta^2\int_{0}^{t}\int\t^\beta(1+s)^{-4}e^{-\frac{1+s}{c_1x^2}}dxds\\[3mm]
		&\di\leq\frac{c_v}{8}\int_{0}^{t}\int\t^\beta\zeta_t^2dxdx +C_0.
	\end{array}
\end{equation}

Combining (\ref{int-t+}
)-(\ref{k5}), (\ref{basic}), (\ref{theta^beta}) and Gronwall inequality leads to 
\begin{equation}
	\begin{array}{ll}
		\di\sup_{0\leq t\leq T}\int(\t^\beta\zeta_x)^2dx + \int_{0}^{T}\int\t^\beta\zeta_t^2dxdt\leq C_0,
		\label{sup-tx}
	\end{array}
\end{equation}
which together with (\ref{cedu}) in particular gives
\begin{equation}
	\begin{array}{ll}
		\di\max_{x\in \mathbb{R}}(\t-2\T)_+&\di\leq\int_{\Omega_2}|\t_x -2\T_x|dx\\[3mm]
		&\di\leq C\int_{\Omega_2}|\zeta_x| + |\T_x|dx \\[3mm]
		&\di \leq C\left(\int_{\Omega_2}\left(\t^\beta\zeta_x\right)^2dx\right)^{1/2} + C\leq C_0.
		\label{the-up}
	\end{array}
\end{equation}
This yields that for all $(x,t)\in \mathbb{R}\times[0,\infty)$,
\begin{equation}
	\di \t(x,t)\leq C_0.
	\label{theta-up}
\end{equation}

Next, multiplying $(\ref{perturb})_3$ by $\zeta^5$ and integrating the resulting equality over $\mathbb{R}$ gives
\begin{equation}
	\begin{array}{ll}
		\di\frac{c_v}{6}\frac{d}{dt}\int\zeta^6dx &\di\leq C\int\frac{\zeta^6(\psi_x+U_x)}{v} +\left|\frac{\zeta^5\phi U_x}{vV}\right|dx + C\int\frac{\zeta^4\zeta_x^2}{v}dx\\[3mm]
		&\di\quad +C\int\left|\frac{\phi\T_x\zeta_x}{vV}\right| dx + C\int\frac{\zeta^5\psi_x^2}{v} + \left|\frac{\phi\zeta^5 U_x^2}{vV}\right|dx\\[3mm]
		&\di\leq C\int(\psi_x^2 + \zeta_x^2)dx +C(M)\int|U_x|dx+ C\int( \phi^2\T_x^2+ \zeta^2U_x^2dx)\\[3mm]
		&\di\leq
		C \int\psi_x^2dx +C(M)\delta+ C\int\zeta_x^2dx,
		\label{dt-int}
	\end{array}
\end{equation}
where in the second inequality we have used 
\begin{equation}
	\begin{array}{ll}
		\di\int\zeta^6dx \leq C\left(\int\zeta^2dx\right)^2\int\zeta_x^2dx \leq C\int\zeta_x^2dx,
	\end{array}
\end{equation}
due to (\ref{basic}) and (\ref{theta-up}). Combining this, (\ref{dt-int}) and (\ref{high-deri}) gives
\begin{equation}
	\di \lim_{t\rightarrow+\infty}\int(\t-\T)^6dx = \lim_{t\rightarrow+\infty}\int\zeta^6 dx = 0.
\end{equation}

Sobolev's inequality shows 
\begin{equation}
	\begin{array}{ll}
		\di\max_{x\in \mathbb{R}}(\t^{\beta+1}-\T^{\beta+1})^2&\di\le C\left(\int(\t^{\beta+1}-\T^{\beta+1})^6dx\right)^{1/4}\left(\int\left|(\t^{\beta+1}-\T^{\beta+1})_x\right|^2dx\right)^{1/4}\\[3mm]
		&\di\le C\left(\int(\t-\T)^6dx\right)^{1/4}\left(\int\left(\t^{2\beta}\zeta_x^2+\T_x^2|\t^\beta-\T^\beta|^2\right)dx\right)^{1/4}\\[3mm]
		&\di\leq C \left(\int(\t-\T)^6dx\right)^{1/4}.
	\end{array}
\end{equation}
Hence there exists some $T_0>0$ such that for all $(x,t)\in \mathbb{R} \times [T_0,+\infty)$, it holds that
\begin{equation}
	\di\t>2^{-\frac{1}{\beta+1}}\T\geq C_0.
	\label{t0-infty}
\end{equation} 

Next, for $p>2$, multiplying (\ref{ener}) by $\t^{-2}(\t^{-1}-2\Theta)^p_+$ and integrating over $\mathbb{R}$. We obtain 
\begin{equation*}
	\begin{array}{ll}
		\di\frac{1}{p+1}\frac{d}{dt}\left(\int(\t^{-1}-2\T)_+^{p+1}dx\right) + \int\frac{u_x^2}{v\theta^2}(\t^{-1}-2\Theta)^p_+dx \\[3mm]
		\di\quad \quad+ \int\frac{\t^{\beta}\t_x^2}{v\t^3}(\t^{-1}-2\Theta)^p_+dx + \int\frac{\t^{\beta}\t_x^2}{v\t^4}(\t^{-1}-2\Theta)^p_+\\[3mm]
		 \di\quad\leq\int\frac{u_x}{v\t}(\t^{-1}-2\Theta)^p_+dx  + \frac{1}{2}\int\frac{\t^{\beta}\t_x^2}{v\t^4}(\t^{-1}-2\Theta)^p_+dx+ C\int(\t^{-1}-2\Theta)^p_+dx\\[3mm]
		\di\quad\leq \frac{1}{2}\int\frac{u_x^2}{v\t^2}(\t^{-1}-2\Theta)^p_+dx + \frac{1}{2}\int\frac{\t^{\beta}\t_x^2}{v\t^4}(\t^{-1}-2\Theta)^p_+dx \\[3mm]
		\di\quad\quad+ C\left(\int(\t^{-1}-2\Theta)_+^{p+1}dx\right)^{\frac{p}{p+1}},
	\end{array}
\end{equation*}
where in the last inequality we have used (\ref{cedu}). Thus, we have 
\begin{equation*}
	\begin{array}{ll}
		\|(\t^{-1}-2\Theta)_+\|^p_{L^{p+1}}\left(\|\t^{-1}-2\T)_+\|_{L^{p+1}}\right)_t\leq C\|(\t^{-1}-2\Theta)_+\|^p_{L^{p+1}}
	\end{array}
\end{equation*}
with C independent of $p$. This in particular implies that there exists some positive constant C independent of $p$ such that 
\begin{equation}
	\di\sup_{0\leq t\leq T}\|(\t^{-1}-2\Theta)_+\|_{L^{p+1}} \leq C_0.
	\label{0-t0}
\end{equation}
Using (\ref{cedu}) and  letting $p\rightarrow\infty$ in (\ref{0-t0}) shows
\begin{equation}
	\di\sup_{0\leq t\leq T}\|(\t^{-1}-2\Theta)_+\|_{L^\infty} \leq C_0,
	\label{0-t0}
\end{equation}
which gives 
\begin{equation}
	\di\theta(x,t)>C_0^{-1}
\end{equation}
for all $(x,t) \in \mathbb{R} \times[0,T_0]$. This combined with (\ref{t0-infty}) and (\ref{theta-up}) finishes the proof of lemma \ref{theta-supbel}.

\
\begin{lemma}
	\label{txxtt}
	There exists a positive constant $C_0$ such that 
	\begin{equation}
		\di \sup_{0\leq t\leq T}\int\zeta_x^2dx + \int_{0}^{T}\int(\zeta_t^2+\zeta_{xx}^2)dxdt\leq C_0.
		\label{zetag}
	\end{equation}
	
\end{lemma}
\textbf{Proof}:\quad First, both (\ref{theta-supbel}) and (\ref{sup-tx}) lead to
\begin{equation}
	\di \sup_{0\leq t\leq T}\int\zeta_x^2dx + \int_{0}^{T}\int\zeta_t^2dxdt \leq C_0.
	\label{t-linfl2}
\end{equation}

Next, it follows from (\ref{ener}) that 
\begin{equation*}
	\frac{(\t^\beta\t_x)_x}{v} = \frac{\t^\beta\t_xv_x }{v}-\frac{u_x^2}{v} + \frac{\t u_x}{v} + \t_t,
\end{equation*}
which together with (\ref{high-deri}), (\ref{theta-fina}), (\ref{t-linfl2}) and lemma 2.1 gives
\begin{equation}
	\begin{array}{ll}
		\di\int_{0}^{t}\int\left|(\t^\beta\zeta_x)_x\right|^2dxds + \int_{0}^{t}\int\left|(\t^\beta\T_x)_x\right|^2dxds\\[3mm]
		\di\quad\leq \frac{1}{2}\int_{0}^{t}\int\left|(\t^\beta\zeta_x)_x\right|^2dxds + C\int_{0}^{t}\max_{x\in \mathbb{R}}(\t^\beta\zeta_x)^2\int\phi_x^2dxds + C\\[3mm]
		\di\quad\leq \frac{1}{2}\int_{0}^{t}\int\left|(\t^\beta\zeta_x)_x\right|^2dxds + C\int_{0}^{t}\max_{x\in \mathbb{R}}(\t^\beta\zeta_x)^2ds + C.
		\label{L3.911}
	\end{array}
\end{equation}

We get by (\ref{theta-fina}) that 
\begin{equation*}
	\begin{array}{ll}
		\di\int_{0}^{t}\max_{x\in \mathbb{R}}(\t^\beta\zeta_x)^2ds&\di\leq \int_{0}^{t}\max_{x\in \mathbb{R}}\left(-2\int_{x}^{\infty}\left(\t^\beta\zeta_x\right)\left(y\right)\left(\t^\beta\zeta_x\right)_x(y)dy\right)ds\\[3mm]
		&\di\leq C\int_{0}^{t}\left(\|\t^\beta\zeta_x\|_{L^2}\|(\t^\beta\zeta_x)_x\|_{L^2}\right)ds\\[3mm]
		&\di\leq C(\varepsilon)\int_{0}^{t}\int(\t^\beta\zeta_x)^2dxds + \varepsilon\int_{0}^{t}\int\left|(\t^\beta\zeta_x)_x\right|^2dxds,
	\end{array}
\end{equation*}
which together with (\ref{L3.911}), (\ref{theta-fina}) and (\ref{high-deri}) implies
\begin{equation}
	\di\int_{0}^{T}\max_{x\in \mathbb{R}}\zeta_x^2dt + \int_{0}^{T}\int\left|(\t^\beta\zeta_x)_x\right|^2dxdt \leq C_0.
	\label{L3.92}
\end{equation}

Finally, since 
\begin{equation*}
	\di \zeta_{xx} = \frac{(\t^\beta\zeta_x)_x}{\t^\beta}-\frac{\beta\t_x\zeta_x}{\t},
\end{equation*}
it follows from (\ref{L3.92}), (\ref{theta-fina}) and (\ref{sup-tx}) that 
\begin{equation*}
	\begin{array}{ll}
		\di\int_{0}^{T}\int\zeta_{xx}^2dxdt &\di\leq C \int_{0}^{T}\int\left|(\t^\beta\zeta_x)_x\right|^2dxdt + C\int_{0}^{T}\max_{x\in \mathbb{R}}\zeta_x^2\int\zeta_x^2dxdt +C\\[3mm]
		&\di\leq C\sup_{0\leq t\leq T}\int\zeta_x^2dx\int_{0}^{T}\max_{x\in \mathbb{R}}\zeta_x^2dt + C\\[3mm]
		&\di\leq C_0,
	\end{array}
\end{equation*}
which together with (\ref{t-linfl2}) gives (\ref{zetag}) and finishes the proof of lemma \ref{txxtt}.

\


\section{Proof of Theorem \ref{theorem2}}
\setcounter{equation}{0}

It is sufficient to show the same a priori estimate as Proposition \ref{prop}. Noticing that
$(V_{\pm}^r,U_{\pm}^r,\T_{\pm}^r)$ satisfies Euler system \eqref{euler} and $(V^{cd},U^{cd},\T^{cd})$
satisfies $\eqref{ns}_1$ \eqref{Theta t}, we rewrite the Cauchy problem \eqref{ns}\eqref{initial} as
\begin{equation}\label{perturb2}
\left\{
\begin{array}{ll}
\di \phi_t-\psi_x=0,\\
\di \psi_t+\left(p-P\right)_x=\tilde{\mu}\left(\frac{u_x}{v}-\frac{U_x}{V}\right)_x+F,\\
\di c_{\nu}\zeta_t+pu_x-PU_x=\tilde{\kappa}\left(\frac{\t^\beta\t_x}{v}-\frac{\T^\beta\Theta_x}{V}\right)_x
+\tilde{\mu}\left(\frac{u_x^2}{v}-\frac{U^2_x}{V}\right)+G,\\
\di(\phi,\psi,\zeta)(\pm\infty,t)=0,\\
\di (\phi,\psi,\zeta)(x,0)=(\phi_0,\psi_0,\zeta_0)(x),\quad x\in \mathbb{R},
\end{array}
\right.
\end{equation}
where
$$
P=\frac{R\T}{V},\quad P_{\pm}=\frac{R\T_{\pm}^r}{V_{\pm}^r},
$$
$$
F=(P_-+P_+-P)_x+\left(\tilde{\mu}\frac{U_x}{V}\right)_x-U_t^{cd},
$$
and
$$
\begin{array}{ll}
G=(p^m-P)U^{cd}_x+(P_--P)(U_-^r)_x+(P_+-P)(U_+^r)_x\\
\di\quad\quad+\tilde{\mu}\frac{U_x^2}{V}+\tilde{\kappa}\left(\frac{\T^\beta\T_x}{V}-\frac{(\T^{cd})^\beta\T_x^{cd}}{V^{cd}}\right)_x\\[3mm]
\di\quad\triangleq G_1 + G_2 + G_3.
\end{array}
$$

Like Lemma \ref{basic-lemma}, we have the following key estimate.
\begin{lemma}\label{basic2}
For $(\phi,\psi,\zeta)\in X([0,T])$, we assume \eqref{same-order} holds, then there exist some positive
constants $C_0$ and $\d_0$ such that if $\d<\d_0$, it follows that for $t\in[0,T]$,
\begin{equation}
\begin{array}{ll}
\di \Big\|\left(\psi,\sqrt{\Phi\left(\frac{v}{V}\right)},\sqrt{\Phi\left(\frac{\t}{\Theta}\right)}\right)(t)\Big\|^2
+\int_0^t\int\left(\frac{\psi_x^2}{\t v}+\frac{\t^\beta\zeta_x^2}{\t^2 v}\right)dxds\\[3mm]
\di+\int_0^t\int P\left(\Phi\left(\frac{\t V}{\T v}\right)+\g\Phi\left(\frac{v}{V}\right)\right)\big((U_-^r)_x+(U_+^r)_x\big)dxds
\di\leq C_0,
\end{array}
\end{equation}
where $C_0$ denotes a constant depending only on $\tilde{\mu}$, $\tilde{\kappa}$, $R$, $c_{\nu}$,
$v_{\pm}$, $u_{\pm}$, $\t_{\pm}$  and $m_0$.
\end{lemma}

\textbf{Proof}: First, multiplying $\eqref{perturb2}_2$ by $\psi$ leads to
\begin{equation}\label{4.17}
\begin{array}{ll}
\di \left(\frac{\psi^2}{2}\right)_t+\left[(p-P)\psi-\tilde{\mu}\left(\frac{u_x}{v}-\frac{U_x}{V}\right)\psi\right]_x
-\frac{R\zeta}{v}\psi_x\\
\di-R\T\left(\frac{1}{v}-\frac{1}{V}\right)\phi_t+\tilde{\mu}\frac{\psi_x^2}{v}
+\tilde{\mu}\left(\frac{1}{v}-\frac{1}{V}\right)U_x\psi_x=F\psi.
\end{array}
\end{equation}

Next, we multiply $\eqref{perturb2}_3$ by $\z\t^{-1}$ to get
\begin{equation}\label{4.18}
\begin{array}{ll}
\di\frac{R}{\g-1}\frac{\z\z_t}{\t}-\left[\tilde{\k}\left(\frac{\t_x}{v}-\frac{\T_x}{V}\right)\frac{\z}{\t}\right]_x
+\frac{R\z}{v}\psi_x+(p-P)U_x\frac{\z}{\t}\\[3mm]
\di+\tilde{\k}\frac{\t^\beta\T\z_x^2}{\t^2v}-\tilde{\k}\frac{\t^\beta\z_x\z\T_x}{\t^2v}-\tilde{\k}\frac{\left(\T^\beta v-\t^\beta V\right)\T_x\T\z_x}{\t^2vV}
+\tilde{\k}\frac{\left(\T^\beta v-\t^\beta V\right)\z\T_x^2}{\t^2vV}\\[2mm]
\di-\tilde{\mu}\frac{\z\psi_x^2}{\t v}-2\tilde{\mu}\frac{\psi_xU_x\z}{\t v}+\tilde{\mu}\frac{\phi\z U_x^2}{\t vV}=G\frac{\z}{\t}.
\end{array}
\end{equation}
Noticing that
\begin{equation}\label{4.19}
\di-R\T\left(\frac{1}{v}-\frac{1}{V}\right)\phi_t=\left[R\T\Phi\left(\frac{v}{V}\right)\right]_t
-R\T_t\Phi\left(\frac{v}{V}\right)+\frac{P\phi^2}{vV}V_t,
\end{equation}
\begin{equation}
\di\left[\T\Phi\left(\frac{\t}{\T}\right)\right]_t=\frac{\z\z_t}{\t}-\T_t\Phi\left(\frac{\T}{\t}\right),
\end{equation}
\begin{equation}
\begin{array}{ll}
-R\T_t&\di=(\g-1)P_-(U_-^r)_x+(\g-1)P_+(U_+^r)_x-p^mU_x^{cd}\\[2mm]
&\di=(\g-1)P(U_-^r)_x+(\g-1)P(U_+^r)_x+(\g-1)(P_--P)(U_-^r)_x\\[2mm]
&\quad\di+(\g-1)(P_+-P)(U_+^r)_x-p^mU_x^{cd},
\end{array}
\end{equation}
and
\begin{equation}
\begin{array}{ll}
\di-R\T_t\Phi\left(\frac{v}{V}\right)+\frac{P\phi^2}{vV}V_t+\frac{R}{\g-1}\T_t\Phi\left(\frac{\T}{\t}\right)
+(p-P)U_x\frac{\z}{\t}\\
\quad=Q_1\big((U_-^r)_x+(U_+^r)_x\big)+Q_2,
\end{array}
\end{equation}
where
\begin{equation}
\begin{array}{ll}
Q_1&\di=(\g-1)P\Phi\left(\frac{v}{V}\right)+\frac{P\phi^2}{vV}-P\Phi\left(\frac{\T}{\t}\right)
+\frac{\z}{\t}(p-P)\\
&\di=P\left(\Phi\left(\frac{\t V}{\T v}\right)+\g\Phi\left(\frac{v}{V}\right)\right),
\end{array}
\end{equation}
\begin{equation}
\begin{array}{ll}\label{q2}
Q_2&\di=U_x^{cd}\left(\frac{P\phi^2}{vV}-p^m\Phi\left(\frac{v}{V}\right)+\frac{p^m}{\g-1}\Phi\left(\frac{\T}{\t}\right)
+\frac{\z}{\t}(p-P)\right)\\[3mm]
&\di\quad+(\g-1)(P_--P)(U_-^r)_x\left(\Phi\left(\frac{v}{V}\right)-\frac{1}{\g-1}\Phi\left(\frac{\T}{\t}\right)\right)\\[3mm]
&\di\quad+(\g-1)(P_+-P)(U_+^r)_x\left(\Phi\left(\frac{v}{V}\right)-\frac{1}{\g-1}\Phi\left(\frac{\T}{\t}\right)\right).
\end{array}
\end{equation}
Adding \eqref{4.19} into \eqref{4.18}, it follows from \eqref{4.19}-\eqref{q2} that
\begin{equation}\label{2-basic}
\begin{array}{ll}
\di\left(\frac{\psi^2}{2}+R\Theta\Phi\left(\frac{v}{V}\right)+\frac{R}{\g-1}\Theta\Phi\left(\frac{\theta}{\Theta}\right)\right)_t
+\frac{\tilde{\mu}\Theta}{\theta v}\psi_x^2+\frac{\tilde{\kappa}\t^\beta \Theta}{\theta^2 v}\zeta_x^2\\
\di+H_{x}+Q_1\big((U_-^r)_x+(U_+^r)_x\big)+\tilde Q=F\psi+G\frac{\zeta}{\theta}
\end{array}
\end{equation}
with $H$ the same as in \eqref{H} and
\begin{equation}\label{tilde-q}
\begin{array}{ll}
\tilde Q&\di=Q_2-\tilde{\k}\frac{\t^\beta\z_x\z\T_x}{\t^2v}-\tilde{\k}\frac{\left(\T^\beta v-\t^\beta V\right)\T_x\T\z_x}{\t^2vV}
+\tilde{\k}\frac{\left(\T^\beta v-\t^\beta V\right)\z\T_x^2}{\t^2vV}\\[3mm]
&\di \quad-\tilde{\mu}\frac{\phi U_x\psi_x}{vV}-2\tilde{\mu}\frac{\psi_xU_x\z}{\t v}+\tilde{\mu}\frac{\phi\z U_x^2}{\t vV}.
\end{array}
\end{equation}
Recalling (iii) in Lemma \ref{rare-pro}, we have
\begin{equation}
\begin{array}{ll}
\di|(P_--P)(U_-^r)_x|\\
\di\leq C\big(|\T^{cd}-\t_-^m|+|\T^{r}_+-\t_+^m|+|V^{cd}-v_-^m|+|V^{r}_+-v_+^m|\big)|(U_-^r)_x|\\
\di\leq C\big(|\T^{cd}-\t_-^m|+|\T^{r}_+-\t_+^m|+|V^{cd}-v_-^m|+|V^{r}_+-v_+^m|\big)\big|_{\Omega_-}
+C|(U_-^r)_x|\big|_{\Omega_c\cap\Omega_+}\\
\di\leq C\d e^{-c_0(|x|+t)},
\end{array}
\end{equation}
which leads to
\begin{equation}
\di|Q_2|\leq C(M)|U_x^{cd}|(\phi^2+\zeta^2)+C(M)\d e^{-c_0(|x|+t)}(\phi^2+\zeta^2)
\end{equation}
and
\begin{equation}
\di|\tilde Q|\leq |Q_2|+\frac{\tilde{\mu}\T\psi_x^2}{4\t v}+\frac{\tilde{\k}\t^\beta\T\zeta_x^2}{4\t^2 v}
+C(M)(\phi^2+\zeta^2)(\T_x^2+U_x^2).
\end{equation}
Note that
\begin{equation}
\begin{array}{ll}
(\phi^2+\zeta^2)(\T_x^2+U_x^2)\leq C(\phi^2+\zeta^2)((\T_x^{cd})^2+(U^r_-)_x^2+(U^r_+)_x^2)\\
\quad \leq C(1+t)^{-1}(\phi^2+\zeta^2)\e^{-\frac{c_1x^2}{1+t}}+C\d Q_1((U^r_-)_x+(U^r_+)_x).
\end{array}
\end{equation}

Since
\begin{equation}
	\begin{array}{ll}
		\label{G3}
		G_3 = \tilde{\kappa}\left(V^{-1}\left((\T_-^r)^\beta\left(\T_-^r\right)_x+(\T_+^r)^\beta\left(\T_+^r\right)_x\right)\right)_x\\[3mm]
		\di\quad\quad + \tilde{\kappa}\left((\T^{cd})^\beta\T^{cd}_x\left(V^{-1}-(V^{cd})^{-1}\right)\right)_x\\[3mm]
		\di\quad\triangleq G_3^1+G_3^2,
	\end{array}
\end{equation}
where 
\begin{equation}
	\label{G31}
	\begin{array}{ll}
		\di |G_3^1| &\di\leq C\left(|(\T_-^r)_{xx}| + |(\T_+^r)_{xx}| + |(\T_-^r)_{x}|^2 +|(\T_+^r)_{x}|^2\right)\\[3mm]
		&\di\quad+C|(\T_+^r)_x|\left(|(V_+^r)_x| + |(V_-^r)_x| +|(V^{cd}_x|\right)\\[3mm]
		&\di\quad+C|(\T_-^r)_x|\left(|(V_+^r)_x| + |(V_-^r)_x| +|(V^{cd}_x|\right),
	\end{array}
\end{equation}
it follows from lemma 2.1 and lemma 2.5 that 
\begin{equation}
	\begin{array}{ll}
		\|G_3^1\|_{L^1}\leq C\delta^{1/8}(1+t)^{-7/8}.
		\label{g31c}
	\end{array}
\end{equation}
And 
\begin{equation}
	\begin{array}{ll}
		\label{g32}
		\di|G_3^2|&\di\leq  C\left(|\T_{xx}^{cd}| + |\T_x^{cd}||V_x^{cd}|\right)\left(|V_+^r-v_+^m| + |V_-^r-v_-^m|\right)\\[3mm]
		&\quad + C|\T_x^{cd}|\left(|(V_-^r)_x + (V_+^r)_x|\right)\\[3mm]
		&\leq C\delta e^{-C_0(|x|+t)}.
		
	\end{array}
\end{equation}
 Following the calculations in \cite{Huang-Li-matsumura}, we obtain
\begin{equation}
\di \|(F,G)\|_{L^1}\leq C\d^{1/8}(1+t)^{-7/8}.
\end{equation}
Therefore, we have
\begin{equation}\label{int-F}
\begin{array}{ll}
\di\int_0^t\int\left(F\psi+G\frac{\z}{\t}\right)dxds\leq C(M)\int_0^t\|(F,G)\|_{L^1}\|(\psi,\zeta)\|_{L^{\infty}}ds\\
\di\quad\leq C(M)\d^{1/8}\int_0^t(1+s)^{-7/8}\|(\psi,\zeta)\|^{1/2}\|(\psi_x,\zeta_x)\|^{1/2}ds\\
\di\quad\leq \int_0^t\int\left(\frac{\tilde{\mu}\T\psi_x^2}{4\t v}+\frac{\tilde{\k}\t^\beta\T\zeta_x^2}{4\t^2v}\right)dxds\\
\di\quad+C(M)\d^{1/6}\int_0^t(1+s)^{-7/6}\left(1+\left\|\left(\psi,\sqrt{\Phi\left(\frac{\t}{\T}\right)}\right)\right\|^2\right)ds.
\end{array}
\end{equation}
Integrating \eqref{2-basic} over $\mathbb{R}\times(0,t)$ and using Gronwall's inequality,
we deduce from \eqref{tilde-q}-\eqref{int-F} that
\begin{equation}\label{7.18}
\begin{array}{ll}
\di \Big\|\left(\psi,\sqrt{\Phi\left(\frac{v}{V}\right)},\sqrt{\Phi\left(\frac{\t}{\Theta}\right)}\right)(t)\Big\|^2
 +\int_0^t\int\left(\frac{\psi_x^2}{\t v}+\frac{\t^\beta\zeta_x^2}{\t^2 v}\right)dxds\\[3mm]
\di+\int_0^t\int Q_1\big((U_-^r)_x+(U_+^r)_x\big)dxds \leq C_0+C(M)\d^{1/6}\\[3mm]
\di+C(M)\d\int_0^t(1+s)^{-1}\int_{\mathbb{R}}(\phi^2+\zeta^2)e^{-\frac{c_1x^2}{1+s}}dxds.
\end{array}
\end{equation}
Finally, due to the fact that
\begin{equation}\label{7.19}
\begin{array}{ll}
\di\int_0^t\int(\phi^2+\psi^2+\zeta^2)w^2dxds
&\di\leq C(M)+C(M)\int_0^t\|(\phi_x,\psi_x,\zeta_x)\|^2ds\\
&\di+C(M)\int_0^t\int(\phi^2+\zeta^2)\big((U_-^r)_x+(U_+^r)_x\big)dxds\\
&\di\leq C(M)+C(M)\int_0^t\int\left(\frac{\t\phi_x^2}{v^3}+\frac{\psi_x^2}{\t v}
+\frac{\t^\beta\zeta_x^2}{\t^2v}\right)dxds\\
&\di+C(M)\int_0^t\int Q_1\big((U_-^r)_x+(U_+^r)_x\big)dxds,
\end{array}
\end{equation}
whose proof can be found in \cite{Huang-Li-matsumura} and $w$ is defined as in Lemma \ref{hlm}, substituting \eqref{7.19} into \eqref{7.18} and choosing $\d$ suitable small imply
\begin{equation}
\begin{array}{ll}
\di \Big\|\left(\psi,\sqrt{\Phi\left(\frac{v}{V}\right)},\sqrt{\Phi\left(\frac{\t}{\Theta}\right)}\right)(t)\Big\|^2
+\int_0^t\int\left(\frac{\psi_x^2}{\t v}+\frac{\t^\beta\zeta_x^2}{\t^2 v}\right)dxds\\[3mm]
\di+\int_0^t\int Q_1\big((U_-^r)_x+(U_+^r)_x\big)dxds
\leq C_0+C(M)\d\int_0^t\int\frac{\t\phi_x^2}{v^3}dxds.
\end{array}
\end{equation}
Thus we can finish the proof of Lemma \ref{basic2} in the similar way as in Lemma \ref{basic-lemma}.
\hfill $\Box$

It is straightforward to verify that the estimates for a single viscous contact wave remain valid for composite waves combining viscous contact waves with rarefaction waves.
Thus, the proof of Proposition \ref{prop} is completed, and Theorem \ref{theorem2} is consequently finished.
          



\end{document}